\documentclass[preprint, 12pt, authoryear]{elsarticle}
\usepackage{amsmath,amssymb}
\usepackage{bbold}
\usepackage{float}
\usepackage{graphicx}
\usepackage{dcolumn}
\usepackage{bm}
\usepackage[colorlinks]{hyperref}
\usepackage{braket}
\usepackage{mathtools}
\usepackage{comment}
\usepackage{xcolor}
\usepackage{physics}
\usepackage{enumerate}
\usepackage{subcaption}

\usepackage{graphicx} 
\usepackage[english]{babel}
\usepackage{physics}
\usepackage{xcolor}
\usepackage{amssymb}
\usepackage{bm}
\usepackage{amsmath}
\usepackage{enumitem} 

\usepackage{caption}
\usepackage{float}
\usepackage{makecell}

\usepackage{tikz}
\usepackage{color}
\usepackage{pgfplots}
\usepackage{varwidth}

\usepackage[T1]{fontenc}
\usepackage[utf8]{inputenc}

\usepgfplotslibrary{fillbetween}
\usetikzlibrary{patterns}
\usetikzlibrary{shapes, arrows}
\usetikzlibrary {arrows.meta}
\usetikzlibrary{calc,math}
\usetikzlibrary{positioning}
\usetikzlibrary{shapes.geometric,backgrounds}

\tikzstyle{scheme step}=[draw, rectangle, rounded corners, align=center, text centered]
\tikzstyle{connector}=[draw, thick, -latex']

\begin{document}
\begin{frontmatter}
\title{Multistart Large Neighborhood Search for the liquefied natural gas transportation and trading over long-term time horizons}

\author{S. Iudin}
\affiliation{Terra Quantum AG, Kornhausstrasse 25, 9000 St. Gallen, Switzerland}

\author{M. Veshchezerova \corref{cor1}}
\affiliation{Terra Quantum AG, Kornhausstrasse 25, 9000 St. Gallen, Switzerland}
\cortext[cor1]{Corresponding author. URL: \url{mv@terraquantum.swiss} (Dr.\ M.\ Veshchezerova)}

\author{K. Tsarova}
\affiliation{Terra Quantum AG, Kornhausstrasse 25, 9000 St. Gallen, Switzerland}

\author{G. Tadumadze}
\affiliation{Terra Quantum AG, Kornhausstrasse 25, 9000 St. Gallen, Switzerland}

\author{V. Shete}
\affiliation{Terra Quantum AG, Kornhausstrasse 25, 9000 St. Gallen, Switzerland}

\author{J.-K. Hao}
\affiliation{LERIA, Universite' d'Angers, 2 Boulevard Lavoisier, 49045 Angers, France}

\author{M. Perelshtein}
\affiliation{Terra Quantum AG, Kornhausstrasse 25, 9000 St. Gallen, Switzerland}

\begin{abstract}

    Liquefied Natural Gas (LNG) transportation is a critical component of the energy industry. It enables the efficient and large-scale movement of natural gas across vast distances by converting it into a liquid form, thereby addressing global demand and connecting suppliers with consumers.
    In this study, we present the Multistart Large Neighborhood Search heuristic for the LNG transportation problem, which involves hundreds of contracts and a planning horizon of two to three years. Our model incorporates several fuel types, LNG sloshing in the tank, and speed- and load-dependent consumption rates. 
    We also consider flexible contracts with LNG volume variability, enabling volume optimizations and multiple discharges. 
A tensor-train optimizer defines the parameters of Mixed Integer Programming (MIP) models, allowing better solution space exploration. 
    On the historic and artificially generated data, our approach outperforms the baseline linear-programming model by $35\%$ and $44\%$, respectively, while the time overhead is only several minutes. 
\end{abstract}


\end{frontmatter}

\section{Introduction}
\paragraph{Background and Motivation}

Maritime transportation is the backbone of international trade; it facilitates the exchange of commodities and raw materials and is essential for industry. Its importance can be attributed to its cost-effectiveness and efficiency in moving large quantities of goods over long distances. Compared to other modes of transport, such as air or land, shipping offers a more economical solution, especially for containerized, bulk, and heavy items\footnote{https://www.seaspace-int.com/sea-vs-air-vs-land-freight-what-is-the-best-method-of-transport-for-you/}. Maritime transportation involves many complex strategic, tactical, and operational decisions: management of terminal operations, route and speed optimizations, fleet deployment, and the overall supply chain management. The increased need for environmental footprint reduction (the key challenge for the sector \citep{UNCTAD2023}) necessitates an increase in energy efficiency.

This paper focuses on \textit{Liquefied Natural Gas Trading and Transportation} problem. 
The topic is particularly relevant today as environmental and geopolitical factors are significantly elevating the role of LNG in the global economy \citep{UNCTAD2023}. Indeed, LNG is a cleaner alternative, producing less carbon dioxide than other fossil fuels for the same energy output\footnote{https://www.eia.gov/}, and LNG spills are less environmentally damaging compared to oil spills \citep{lehr2017comparative}. Furthermore, LNG contributes to the diversification of energy supply sources and decreases reliance on gas from any single country or region, offering an advantage over pipeline-dependent gas.

\paragraph{Problem overview} 
In the traditional ship scheduling problem \citep{Appelgren1969}, the vessel fleet transports fixed or optional cargoes from corresponding sources to destinations while respecting vessel capacity constraints and time windows for pickup and delivery. Our LNG transportation problem differs from this traditional formulation in several ways: we allow for a \textit{heterogeneous vessel fleet} with possible restrictions on vessel availability, the fuel consumption depends on payload and the speed \citep{li_discretized_speeds}, vessels may use the transported LNG as fuel, among others. Crucially, we don't transport cargoes in the traditional sense but rather select from potential trading opportunities — each tied to a specific location and time window — to maximize profit. This refined problem formulation with detailed vessel characteristics more accurately mirrors real-world conditions \citep{psaraftis_horse} and aligns with the operational challenges faced by the energy providers across the globe.

In our trading model, we can distribute a single load across multiple deliveries (fulfilling smaller island contracts) to enhance profitability. The major complexity comes from the restriction on volumes onboard during the navigation: when the vessel sails, the LNG tank should be either almost full (at least $75\%$ of maximum capacity) or nearly empty (no more than $25\%$ of $V_{\max}$) to prevent \textit{sloshing}: a liquid movement inside a partially loaded container that may challenge vessel's stability. 

\paragraph{Contributions}
 
In this work, we consider a real-world LNG transportation case. To our knowledge, the combination of technical features is unique in scientific literature. The technical peculiarities include a long-range planning horizon (between two and three years), LNG sloshing, multiple discharges, fuel consumption dependency on speed and payload, and other details. To deal with the high complexity of this problem, we propose an effective heuristic approach, which relies on three complementary components: black-box optimization, using \textit{tensor-train-based solver} \citep{zheltkov2020global}, and two consecutive arc-flow-based models.

On the production data, our method yields a 35\% profit increase compared to the baseline linear programming solution that does not allow distributing a single load into multiple deliveries. Crucially, we do so without a significant runtime overhead—our algorithm takes only several minutes more than the baseline. 

We believe that our work, introducing a novel technique involving tensor-train-based optimization for a real-world case study, is significant for maritime transportation and optimization communities. 

\paragraph{Structure of the paper} In Section \ref{sec:lit_review}, we review the relevant literature on LNG transportation.
Section \ref{sec::problem_description} describes the problem, including the input structure, taken assumptions, and optimization objective. In Section \ref{sec:solution_approach}, we introduce our solution approach: the two components of the Large Neighborhood Search are presented in Sections \ref{sec:big_pair} and \ref{sec:small_cargo_model}, and the initial point selection strategy is detailed in Section \ref{sec:black-box}. The MIP formulation for the problem is presented in \ref{sec:mip}. In Section \ref{sec:results}, we present performance analysis for our algorithm on artificial data generated based on a real commercial case. Section \ref{sec:conclusion} summarizes the contributions of this study and discusses possible future steps.

\section{Related works}
\label{sec:lit_review}

A growing research activity in LNG transportation problems debuted in late 2000, with the first study \citep{gronhaug2009_supply_chain} introducing the \textit{Inventory Routing Problem} for LNG, called \textit{LNG-IRP}. Since then, the number of publications on the topic has exploded: both the decimal Maritime Transportation review from 2012 \citep{christiansen2013_review} and the recent review on the Maritime Inventory Routing (MIRP) \citep{review_2023} devote special sections for the LNG. In the last decade, LNG became the most frequent commodity in MIRP studies \citep{review_2023}.

The LNG supply chain consists of multiple stages. In the \textit{upstream} phase, \textit{liquefaction plants} produce LNG, which can be stored in tanks before delivery. The \textit{midstream} phase focuses on the transportation of LNG. Finally, the \textit{downstream} phase covers the regasification of LNG and its further distribution. Regasification plants might also store the LNG before converting it back into gas form. While some works consider all parts of the chain \citep{gronhaug2009_supply_chain, gronhaug_bnp, LNGScheduler}, most of the research is focused on the upstream and midstream sections \citep{review_2023}. In the latter case, the deliveries are modeled as cargoes with specified volume ranges and time windows (e.g., \citep{Halvorsen_ADP}). Usually, the delivery is accomplished with the company-owned fleet, but some models allow chartering \citep{rakke_delivery_patterns, Rakke_ADP_and_rolling_horizon}. 

A wide range of models exists for LNG-IRP, varying based on their decision-making flexibility and closeness to reality.

\subsection{Problem variants and Extensions}

\paragraph{Planning horizon} In LNG supply chain, the optimization problems appear in \textit{strategic}, \textit{tactical} and \textit{operational} time-scales. 
In our problem, the operation horizon spans \textit{several years}, which brings it close to the \textit{Annual Delivery Program (ADP)}, where one wants to plan production and delivery for long-term contracts over the year. ADP belongs to the tactical level but can be used to model “what-if” scenarios for strategic decisions \citep{LNS_for_LNG-IRP, papaleonidas_decision_support_tool} such as the fleet size and composition \citep{bittante_small_scale}. 

The pioneering work in ADP \citep{Rakke_ADP_and_rolling_horizon} introduced the rolling horizon heuristic that solves the MIP formulation over short time periods. Since then, the problem has been addressed in multiple works \citep{Multu_multiple_discharges, Halvorsen_ADP, Haug_ADP, rakke_delivery_patterns, stalhane_adp_heuristic} suggesting different models and solution approaches. The work \citep{fix-and-relax_heuristic_LNG-IRP} solves the problem on a half-year planning horizon. 

The observed complexity of MIP-based approaches for the LNG delivery problem, where decisions are indexed with discrete time steps, scales badly with the time horizon \citep{rakke_delivery_patterns, Haug_ADP}. 
However, they can still be used for annual planning if integrated inside the rolling-horizon heuristic as the algorithm to solve subproblems. 

Contrary to short-term plans, a realistic annual program should reserve some time for \textit{vessel maintenance}. The maintenance requirement was considered in \citep{Haug_ADP, Multu_multiple_discharges, msakni_short_term, LNGScheduler}. We don't explicitly manage maintenance, but it could be integrated into our model quite straightforwardly. 

\paragraph{Boil-off}
The peculiarity of LNG cargo is that, contrary to other goods, a fraction of LNG evaporates every day. This phenomenon is called \textit{boil-off}, and the works \citep{decomposition_LNG-IRP, fix-and-relax_heuristic_LNG-IRP, bittante_small_scale} highlighted the complexity it induces. In some modern vessels, the LNG boil-off effect may be partially reduced with re-liquefaction; this option is considered in \citep{qatar_case_study}; some vessels in our case study also have this capability. Moreover, we consider the boil-off rate variability on load (laden loads typically lead to higher boil-offs); to our knowledge, the only other work considering this variability is \citep{LNGScheduler}.

To mediate the boil-off, some works \citep{gronhaug2009_supply_chain, decomposition_LNG-IRP, qatar_case_study} suggest using the cargo as fuel. We also integrate this option, and our vessels can navigate either fully on LNG, or on the usual fuel, or on a combination of both. 

\paragraph{Multiple discharges}

In most of the works, \citep{Rakke_ADP_and_rolling_horizon, stalhane_adp_heuristic, LNS_for_LNG-IRP, Halvorsen_ADP, papaleonidas_decision_support_tool, Haug_ADP}, LNG carriers operate in a shipload fashion: they are fully loaded in the liquefaction plant and fully unloaded in the regasification plant. 

While only one work allows multiple loads \citep{Haug_ADP}, multiple discharges are subject to active research. In principle, allowing partial unloads increases the options for contract satisfaction. However, the economic viability of multiple discharges remains an object of debate \citep{review_2023}. So does the technical feasibility: when a vessel transports liquids, the movement of liquid in a partially loaded tank, called \textit{sloshing}, may challenge its stability \citep{sloshing_calcium_carbonate, Multu_multiple_discharges}. For this reason, the models in \citep{decomposition_LNG-IRP, LNS_for_LNG-IRP, gronhaug2009_supply_chain} allow unloading only full tanks; multiple discharges are possible only when the LNG carrier has multiple tanks. The studies \citep{gronhaug2009_supply_chain, gronhaug_bnp, decomposition_LNG-IRP, qatar_case_study} report that a maximum of two unloads is viable.

The works \citep{bittante_small_scale} on small LNG carriers serving archipelagos and \citep{LNGScheduler} on a real-world LNG supply chain consider multiple discharges from a single-tank but \textit{ignore} sloshing. We also operate on single-tank carries, but contrary to the previous works, we explicitly address the sloshing effect by introducing \textit{forbidden volume zones}. Indeed, the sloshing is dangerous only if the volume onboard is within some range. To ensure that the carried volume never enters \textit{the forbidden volume zone}, our model constantly tracks the carried volume on each vessel, which raises the complexity of the underlying optimization problem. To our knowledge, our study is the first to consider partial tank unloads \textit{and} sloshing together. 

In \citep{qatar_case_study}, the possibility of multiple discharges implied the reduction of operational costs by $23.4\%$ for a one-month plan. In contrast, in the annual programs considered in \citep{Multu_multiple_discharges}, average cost savings are around $4.2\%$. In our case, allowing multiple discharges leads to an impressive $35\%$ profit increase on production data provided by the anonymous industrial partner.

\paragraph{Spot contracts} Natural gas production inherently fluctuates, and its consumption can also vary based on specific demands. Due to this variability, producers often engage in spot sales to offload surplus LNG. Trading on the spot market can have a more significant impact on profit margins compared to transportation costs minimization \citep{LNGScheduler}.

While spot contracts are particularly relevant in the short-term optimization, \citep{msakni_short_term}, the works on the Annual Delivery \citep{LNGScheduler, Rakke_ADP_and_rolling_horizon, rakke_delivery_patterns, fix-and-relax_heuristic_LNG-IRP} explicitly integrate spot contracts in their LNG-IRP models. The work \citep{papaleonidas_decision_support_tool} goes further by considering exclusively spot contracts in a \textit{tramp shipping} manner: in their model, the fleet is used to optimize profit by taking pickup and delivery requests. 

Most of the existing studies consider cost minimization, leaving spot contracts as a unique option for extra revenue generation. In our work, we follow a different approach where we aim to maximize the total profit, somewhat similar to the approach taken in \citep{LNGScheduler}. 
As in our model, in \citep{LNGScheduler} LNG prices vary in time and location, leaving opportunities for higher profit margins.

\paragraph{Speed optimization}

While speed optimization is more and more present in the research on maritime transportation \citep{christiansen2013_review, li_discretized_speeds}, most of the works on LNG 
\citep{decomposition_LNG-IRP, fix-and-relax_heuristic_LNG-IRP, bittante_small_scale} consider constant speeds and do not include speed-dependent fuel costs. 

Speed optimization may reduce calls to charter vessels and allow for the provision of more spot contracts. Besides, it may reduce fuel consumption and carbon footprint by avoiding unnecessary rapid sails or replacing sails that are too long \citep{qatar_case_study}. For these reasons, the work by \citet{msakni_short_term} studies the model, which allows variable speeds; the recent paper \citep{Haug_ADP} includes the speed-dependent fuel costs in the model, and in \citep{papaleonidas_decision_support_tool}, the fuel costs are subtracted at the end. The work \citep{Haug_ADP} demonstrated that with speed optimization, the vessel uses the entire speed spectrum, highlighting the limitations of the fixed-speed approach.

In our study, we consider not only the dependency of fuel consumption on speed but also on payload, which is crucial for accurate modelling \citep{psaraftis_horse}, especially when LNG cargo is utilized as fuel. Previous research has often overlooked this aspect.

\subsection{Solution techniques}

The LNG-IRP problem is NP-hard \citep{coelho2014thirty}; the complexity on real-world instances heavily depends on the specific model used.

\paragraph{Exact methods}
The original work on LNG-IRP \citep{gronhaug2009_supply_chain} presented two Mixed Integer Linear Program (MILP) formulations - the so-called \textit{arc-flow} and \textit{path-flow} models. In a nutshell, variables in arc-flow models correspond to one sailing leg (a vessel $v$ travels from port $i$ to port $j$, starting on day $d$). In contrast, in the path-flow models, variables correspond to entire journeys of the vessels over the whole time horizon. The path-flow formulations usually have too many variables for an explicit enumeration but are well-suited for column-generation-based approaches \citep{Appelgren1969, gronhaug_bnp}.

Both works \citep{LNGScheduler, Haug_ADP} use arc-flow formulations for their rich models. Arc-flow models scale badly with the timescale: while on test instances over 6 months, the MIP in \citep{Haug_ADP} is solved to $1\%$ gap in one hour, no solution was found for the largest instance with a timespan exceeding a year. 

For the pickup and delivery variation, the work \citep{papaleonidas_decision_support_tool} introduces an arc-flow MILP model. The MILP uses variables for ballast and laden trips, similar to the baseline model in our work. 
The Branch \& Price algorithm for the path-flow formulation of the maritime pickup and delivery problem is presented in \citep{Homsi2020_UHGS}.

\paragraph{Decompositions}

Contrary to the Branch \& Price approaches, where the subproblem solution generates vessel paths for the entire planning period, the works \citep{rakke_delivery_patterns, decomposition_LNG-IRP} use variables for small fractions of the whole journeys, usually involving several legs. The work \citep{rakke_delivery_patterns} splits the routes into \textit{delivery patterns}, and the presented Branch \& Price \& Cut method on this decomposition outperforms the heuristics from \citep{stalhane_adp_heuristic, Rakke_ADP_and_rolling_horizon}. The work \citep{decomposition_LNG-IRP} decomposes the schedule on pre-generated \textit{duties} comprising one load and two deliveries, further improving over \citep{rakke_delivery_patterns} on the solution time and quality. 

A slightly different decomposition angle is taken in \citep{Halvorsen_ADP}, where one first decides on assignment and sequencing of deliveries and then optimizes the precise delivery moments.

\paragraph{Heuristics}

The first work on the Annual Delivery Planning \citep{Rakke_ADP_and_rolling_horizon} introduces a Rolling-Horizon heuristic that separates the planning period into small chunks and sequentially solves the corresponding reduced MIP formulation. A multi-start local search heuristic from \citep{stalhane_adp_heuristic} achieves the same solution quality but in a shorter computational time. Similar to \citep{Rakke_ADP_and_rolling_horizon}, the fix-and-relax heuristic \citep{fix-and-relax_heuristic_LNG-IRP} decomposes the problem over time, while the large-neighborhood search heuristic \citep{LNS_for_LNG-IRP} also fixes some variables in the solution and solves the MIP formulation on the remaining variables. Interestingly, in a variety of ways, several works  \citep{Rakke_ADP_and_rolling_horizon, fix-and-relax_heuristic_LNG-IRP, LNS_for_LNG-IRP, msakni_short_term} use MIP solvers to address subproblems. 

Contrary to the traditional decomposition on time, in our two-stage approach, we prefix big trades on the first stage and insert small discharges in between by solving another MILP. Our algorithm is also a multistart large-neighborhood search with a problem-specific neighborhood that includes routes with small discharges. The work \citep{Multu_multiple_discharges} solves the problem variant with multiple discharges with a custom vessel-routing heuristic. 

Some research works explored powerful metaheuristic approaches for the LNG shipping problem. The work \citep{qatar_case_study} uses the Particle Swarm optimization on a model with environmental impact. The work \citep{Homsi2020_UHGS} on general tramp shipping introduces the state-of-the-art memetic algorithm for the pickup and delivery variation, as the one considered in \citep{papaleonidas_decision_support_tool}. For the same problem, the work \citep{harwood_quantum} presents a QUBO formulation, enabling the application of such quantum algorithms as Quantum Adiabatic Evolution \citep{farhi_qa}, Quantum Annealing \citep{Dwave_issues}, and Quantum Approximate Optimization Algorithm \citep{QAOA_original}.

In contrast to the existing works, our optimization model explicitly addresses sloshing. This makes our model more complex than usual: we need to track the volumes of LNG onboard with extra variables and constraints. Besides, our model considers multiple discharges, vastly increasing the number of feasible visit sequences for each vessel.
We utilize two optimization models in sequence to address this complexity. The first model ignores small contracts and considers only sequential loads and deliveries ("big" enough to avoid sloshing). This model provides a starting point for a neighborhood search, where the neighborhood consists of all plans where additional small contracts can be visited between big buy and sell contracts found at the previous stage.  Thus, we track volumes only in the second stage, moving the most complex part of the problem to the smaller subproblem. The black-box optimizer diversifies the starting points: it suggests different parameters for the first model with the objective of maximizing the overall two-stage profit. In our context, "multistart large neighborhood search" refers to exploring complex neighborhoods via sequential MILP subproblems rather than the classical destroy‑and‑repair scheme, while still following the core principle of optimizing over large structured neighborhoods.

\section{LNG shipping problem statement}
\label{sec::problem_description}
In this section, we formalize the LNG shipping problem as an optimization problem and discuss the assumptions and limitations.

\subsection{Input description}
In our LNG trade and transportation problem, we deal with a heterogeneous fleet of vessels, where each vessel has the following features:

\begin{itemize}
    \item \textit{Available capacity} of the LNG storage, which may not be exceeded
    \item \textit{Forbidden volume zone}, usually between $25\%$ and $75\%$ of the vessel gross capacity. For safety reasons, the amount of LNG onboard must never enter this zone while sailing.
    \item \textit{Fuel consumption table}, specifying fuel consumption for each feasible combination of fuel type, vessel speeds, and payload. Generally, the vessel may move either on traditional fuel or on LNG. Some vessels may also use both LNG and fuel in a combined mode.
    \item \textit{Rent interval}, when the vessel is available for voyages
    \item \textit{Initial and final locations}, where the vessel should start and end its journey.
\end{itemize}

We consider only optional \textit{sell} and \textit{buy} contracts. Thus, we freely choose both consumers and suppliers from the list of available contracts, each having the following attributes:

\begin{itemize}
    \item \textit{Time window} indicating the period when the contract is accessible. In the production data, most contracts are available for only one day, but some have a time window duration of up to a month.
    \item \textit{Volume boundaries} indicating the minimum and maximum acceptable trading volume of LNG
    \item \textit{Port} specifying the geographical location
\end{itemize}

The price of each contract depends on the respective port and changes over time.

\subsection{Assumptions}

We made the following assumptions:

\begin{itemize}
    \item Each loading or unloading operation takes exactly $24$ hours, independently of the amount of LNG onboard
    \item The vessel moves with fixed speed and fuel mode on the path from one cargo to another until it reaches the destination port, and then it idles with only LNG boil-off wastes, waiting for the loading/discharge time
    \item Contracts may not be split: each contract may be taken strictly by one vessel
    \item Unit prices for all contracts are known in advance
    \item Each vessel has exactly one LNG tank, for which we must avoid sloshing
\end{itemize}

All of these assumptions are widely accepted in the literature, and we consider them close to reality, except for the single tank assumption. Usually, LNG tankers have 4-6 independent tanks, but for simplicity, we follow the approach taken in \citep{LNGScheduler} and consider the case with one tank.

\subsection{Decisions and objective}

We aim to assign optional buy and sell contracts to ships and specify i) a trading volume for each selected contract and ii) the time at which we process it. In addition, we have to decide the fuel type used to navigate between each pair of visited ports. We do not consider the rental cost and don't allow chartering, so our cost function doesn't include fleet deployment expenses. The goal is to build operationally feasible fleet routes that maximize the company's profit. 

Although our problem has many technical peculiarities, such as different fuel types and speed optimization, the most complex part is LNG volume management. To remain feasible, the amount of LNG onboard must never exceed the capacity or enter the forbidden volume zone. Our LNG shipping problem combines this requirement with the possibility of multiple discharges, which makes continuous volume variables necessary for the modelling.

A Mixed Integer Program (MIP) formulation (with simpler assumptions) in \ref{sec:mip} clearly demonstrates how the volume management makes the problem non-linear, so that for a long time horizon it is intractable unless decomposition is used. 
In the MIP formulation, we use binary variables $x_{i,k,q}$ taking value 1 if vessel $k$ serves contract $i$ at the $q$-th position in its route, and 0 otherwise. We restrict the solution space by fixing the visit times of contracts to the midpoint of their respective time windows, and the vessels can sail only on LNG. However, even after the simplifications, this model remains too challenging on a several-year time horizon, and the commercial solver on our machine fails due to the lack of memory. Therefore, we split the time horizon into half-year time periods.
After the decomposition, on the production data, \textit{Gurobi}\footnote{https://www.gurobi.com/} finds the sub-optimal solutions for the resulting subproblems in 45 minutes.

\section{Solution approach}\label{sec:solution_approach}

Our principal contribution is a three-component approach to the problem, where two optimization models are solved sequentially -- the second model searches the large neighborhood around the solution of the first model-- and a tensor-train-based black-box optimizer modifies the parameters of the first model to explore different regions of the solution space. Our algorithm takes less than 7 minutes to optimize the production instance using an open-source solver SCIP \citep{SCIP, SCIP_Optimization_Suite} for subproblems. Section \ref{sec:big_pair} describes the \textit{Big-pairs selection model} - without multiple discharges - which provides an initial point for the Large Neighborhood Search. Section \ref{sec:small_cargo_model} explains the \textit{Small contracts insertion model}, which adds small-volume contracts between large purchases and sales. 

\subsection{Big-pairs selection model}\label{sec:big_pair}

The big-pairs selection model relies inherently on a simple arc-flow model used for ship scheduling in \citep{Appelgren1969}. 

\subsubsection{Preprocessing and Sets}

In this work, a \textit{trip} denotes a tuple consisting of a travel arc with some extra attributes.

For our first model, we generate trips of 4 types:
\begin{itemize}
    \item \textbf{Initial Trip} -- the trip which connects the initial vessel location to a loading cargo
    \item \textbf{Laden Trip} -- the trip which is completed while moving from supplier to consumer
    \item \textbf{Ballast Trip} -- the trip which is completed by a vessel while moving from a consumer to a new supplier
    \item \textbf{Final Trip} -- the trip which connects a vessel discharge to the final vessel location
\end{itemize}

All the obtained trips for all the vessels $v$ are indexed by the index $t\in \mathcal{T}$, where $\mathcal{T}$ is the set of all generated trips. Each trip has a specific start contract with a port and date ($p_{start}^t$ and $d_{start}^t$), an end contract with a port and date ($p_{end}^t$ and $d_{end}^t$), a vessel $v_t$ and a vessel speed $s_t$, a fuel type and a traded LNG volume. The trip $t$ generates profit $P_t$. We compute $P_t$ by calculating the optimal sold volume for the given pair of contracts and dates. For ballast trips, $P_t$ represents fuel consumption. 

We generate a complete set of feasible trips that satisfy all the volume and time restrictions (including sloshing condition), vessels' characteristics and limitations, consumption dependency on load, etc. For each laden trip, we discharge the maximum possible amount of LNG, maximizing the local profit from the trade. However, we allow laden trips with minimum buy volume larger than maximum sell volume (minus boil-off and fuel consumption), as the extra volume may be distributed across small-volume contracts in the Small-contracts insertion phase.

In the following, we use the following notations:
\begin{itemize}
    \item $\mathcal{O}$ -- set of all contracts: buy and sell
    \item $\mathcal{V}$ -- set of all vessels
    \item $\mathcal{C}_c$ -- set of all laden trips, either leaving from the contract $c$ or finishing in this contract.
    \item $\mathcal{M}_c$ -- set of all time steps where the contract $c$ can be served
    \item $\mathcal{A}_{v,c,d}$ -- set of all vessel $v$ trips ending at the contract $c$ on day $d\in\mathcal{M}_c$
    \item $\mathcal{D}_{v,c,d}$ -- set of all vessel $v$ trips departing from the contract $c$ after serving it on day $d$.
    \item $\mathcal{I}_v$ -- set of initial trips for the vessel $v\in\mathcal{V}$
    \item $\mathcal{F}_v$ -- set of final trips for the vessel $v\in\mathcal{V}$
\end{itemize}

The sets $\mathcal{O}$ and $\mathcal{M}_c$ are directly provided in the input data, and we compute other sets in the preprocessing step.

\subsubsection{Penalized costs}\label{sec:penalized}

The standard arc-flow model from \citep{Appelgren1969} does not account for the possibility of small (less than $50\%$ of the volume) discharges that may change the preferences between contract pairing in the optimal solution. To mitigate this unawareness, we alter the profits $P_t$ with a couple of trip-specific penalties.

For each laden trip, we introduce the number $w_t$, which corresponds to the over-delivery at the selling point that would happen if no small discharges appeared between the load and the discharge. The value $w_t$ is equal to the difference between the minimum buy volume and the maximum sell volume minus the volume reduction from the sailing. For both laden and ballast trips, we additionally introduce an integer metric $y_t$, which represents the potential for a small discharge between the trip's endpoints. We evaluate $y_t$ by considering how many available sell contracts the vessel $v_t$ can visit between $d_{start}^t$ and $d_{end}^t$ while respecting operational constraints. The metrics $w_t$ and $y_t$ are integrated into the trip's revenue (in USD) $P_t$ with penalty $O$ for the over-delivery in the big-pair solution and the reward $R$ for the multi-destination potential. Thus, the value of the trip $t$ is:

\begin{equation}
\label{eq:penalized_cost}
\pi_t=P_t-O w_t + R y_t.
\end{equation}

The values of $O$ and $R$ impact the amount of LNG overdelivery in the modified Big-pair model solution. We want to allow just enough LNG overdelivery to sell it in small quantities, but not too much: in the final solution, overdelivery is forbidden.

\subsubsection{Variables}

\begin{itemize}
    \item $X_{t}$ -- binary variable indicating if the trip $t\in \mathcal{T}$ is taken in solution, or not:
    \begin{equation}
        X_t = \begin{cases}
        1 & \text{trip $t$ is taken}\\
        0 & \text{trip $t$ doesn't appear in the solution.}
        \end{cases}
    \end{equation}
\end{itemize}

\subsubsection{ILP formulation}

\begin{align}     
&\hspace{5em} \sum_{t \in \mathcal{T}}\pi_{t}X_{t} \to \max &\label{cost}\\ 
&\hspace{5em} \sum_{t \in \mathcal{C}_{c}}X_{t} \leq 1 & \forall c \in \mathcal{O} \label{only_once}\\
&\hspace{5em} \sum_{t \in \mathcal{A}_{v,c,d}}X_{t} = \sum_{t \in \mathcal{D}_{v,c,d}}X_{t} & \forall v \in \mathcal{V},
\forall c \in \mathcal{O}, \forall d \in \mathcal{M}_c \label{arrive_departure}\\
&\hspace{5em} \sum_{t \in \mathcal{I}_v} X_{t} \leq \sum_{t \in \mathcal{F}_v} X_{t} \leq 1 & \forall v \in \mathcal{V}
\label{edge}\\
&\hspace{5em} X_t\in \{0, 1\} & \forall t \in \mathcal{T}
\end{align}

Equation (\ref{cost}) represents the profit for the company from LNG trading. Constraint (\ref{only_once}) ensures that each contract is served not more than once. The flow constraint (\ref{arrive_departure}) guarantees the continuity of the vessel movement. Edge constraints (\ref{edge}) assert each vessel is not used at all or leaves the initial port one time and arrives in the final port. Together with the flow preservation constraint (\ref{arrive_departure}), it also guarantees that the vessel follows a specific route and can't be in multiple locations at the same time.

\subsection{Small discharges insertion model}\label{sec:small_cargo_model}

After solving the big-pair model, we try to insert additional discharges from the remaining sell contracts between already assigned big trades. We also optimize the volumes for the "big" contracts, considering them as continuous variables. In this second stage, we search a large neighborhood around our "big-pair" solution, where neighbors differ from the solution in small discharges.

\subsubsection{Preprocessing and Parameters}

From the big-pairs model solution, we extract the following sets:

\begin{itemize}
    \item $\mathcal{P}_v = (b^v_0, s^v_0), (s^v_0, b^v_1), ...$ -- the ordered sequence of big laden and ballast trips assigned to the vessel $v$. The symbols $b^v_i$ ($s^v_i$) denote buy (sell) contracts taken by the vessel. The number of pairs $\abs{P_v}$ for a vessel $v$ equals the number of trips travelled by the vessel $\abs{T_v}$ minus the initial and the final trip.
    \item $B_v$ -- a set of buy contracts in $\mathcal{P}_v$
    \item $S_v$ -- a set of sell contracts in $\mathcal{P}_v$
    \item $\mathcal{B}$ - all contracts visited in the big-pair model solution
\end{itemize}

We eliminate the already assigned contracts from the data; we filter the remaining sell contracts and keep only the ones allowing small volumes in a new set $\mathcal{S}$. For simplicity, we assume that the contracts from $\mathcal{S}$ may be visited only in the middle of their time windows. For each trip $t$ in the big-pair solution, we filter the set $\mathcal{S}$ to get small-volume discharges that can be served between the endpoints of $t$ by the vessel $v_t$. We call the set of these small contracts $C_t^v$.

For each trip $t$ in the big-pair solution, we compute arcs that may appear in alternative sailing routes relaying the starting and the ending locations of the trip. These involve arcs leaving the start point of $t$ towards contracts in $C_t^v$, the arcs between small deliveries in $C_t^v$, and the arcs leaving the contracts in $C_t^v$ towards the end contract of the trip $t$. The resulting arcs (legs), together with the arc relaying trip's endpoint, form a set of sailing legs that we call $L_t^v$. The LNG consumption of a leg $l\in L_t^v$ is denoted by $lng_l$. We assume that the sailing legs in $L_t^v$ use the same fuel as in the trip $t$. The vessel's speed is the minimal speed from the consumption table for the chosen propulsion mode that allows sailing between the leg endpoints in time. 

In what follows, we'll use the notations:

\begin{itemize}
    \item $vol_v$ - volume of the vessel $v\in \mathcal{V}$
    \item $C^v$ - all small contracts from $\mathcal{S}$ that can be visited by the vessel $v$. It's the union of all sets $C^v_t$ for all trips assigned to the vessel $v$ in the "big-pair" solution
    \item $L^v$ - all legs that can be sailed by the vessel $v$ in a solution with small discharges in the vicinity of the "big-pair" solution. It's the union of the sets $L^v_t$, plus legs for the initial and the final trips. 
\end{itemize}

\subsubsection{Variables}

For every big (buy or sell) contract $c \in \mathcal{B}$ we introduce two volume variables: 
\begin{itemize}
    \item $V_c \in [V_c^{\min}, V_c^{\max}]$ -- the volume traded through in the contract
    \item $F_c$ -- the out-of-contract LNG purchase in the port used for sailing
\end{itemize}

For each big sell $c \in \mathcal{S}$, we also introduce the variable $W_c$ for the over-delivery volume.

 The cost coefficient for $V_c$ is the unit volume price $u_c$ for the contract $c$ (multiplied by $-1$ for the supply contracts), for $F_c$ - the LNG price in the port at the moment when the contract $c$ is visited (we'll refer to it as free price, and it may differ from the price of the contract due to contract peculiarities). For $W_c$, we introduce an over-delivery penalty  - $\omega$ - we pick it high enough to discourage LNG over-delivery in all cases where a feasible solution without LNG waste exists. For each pair $p$ in $\mathcal{P}_v$ we add $2$ auxiliary variables:
 \begin{itemize}
     \item $V_p^{start}$ -- the partial volume at the start of the pair $p\in\mathcal{P}_v$
     \item $V_p^{end}$ -- the partial volume at the end of the pair $p\in\mathcal{P}_v$
 \end{itemize}

For each small discharge $c$, we introduce one continuous volume variable for each vessel $v$ that can visit it -- $ V^v_c$. We force $V^v_c$ to zero if the vessel $v$ doesn't visit the contract with a linear constraint. 

Finally, for each sailing leg in the union of $L^v$, we add one binary variable $l$, indicating that the leg appears in the final routes.

\subsubsection{Mixed Integer Linear Program}

The objective is to maximize the revenue:

\begin{align}
    \sum_{v\in\mathcal{V}}\left(-\sum_{c\in S_v}(u_cV_c + \omega W_c)  + \sum_{c\in B_v}u_cV_c - \sum_{c\in B_v\cup S_v} f_cF_c+\sum_{c \in \mathcal{S}} \sum_{v: c \in C^v} u_cV^v_c \right)\label{small_discharges_objective}
\end{align}

The flow preservation constraints impose that exactly one leg enters and leaves each preselected big contract:

\begin{align}
& \sum_{(c, c')\in L^v} l_{(c, c')}  = \sum_{(c', c)\in L^v} l_{(c', c)} =1, & \forall v\in\mathcal{V}, \quad \forall c \in  B_v \cup S_v \label{flow_preservation1}
\end{align}

The small discharges are visited at most once, and if a vessel enters the contract $c \in \mathcal{S}$, it has to leave it. 

\begin{align}
& \sum_{(c, c')\in L^v} l_{(c, c')} = \sum_{(c', c)\in L^v} l_{(c', c)}, & \forall v \in \mathcal{V}, \quad \forall c \in  C^v \label{enter_exit_legs_equality}
\end{align}

\begin{align}
    &\sum_v \sum_{(c, c')\in L^v} l_{(c, c')} \leq 1 & \forall c \in \mathcal{S} \label{small_only_once1}\\
    &\sum_v \sum_{(c', c)\in L^v} l_{(c', c)} \leq 1 & \forall c \in \mathcal{S}\label{small_only_once2}
\end{align}

If the vessel visits the small discharge contract $c$, the contract's traded volume should be within the contract's bounds; otherwise, the volume is zero. 
\begin{align}
& V_c^{\min}\sum_{(c, c')\in L^v} l_{(c, c')} \leq V_c^v\leq V_c^{\max}\sum_{(c, c')\in L^v} l_{(c, c')}, & \forall v \in \mathcal{V},\quad \forall c \in  C^v \label{if_not_visited_then_zero}
\end{align}

The volumes of LNG onboard change due to purchases and sales, and also decrease during sailing due to LNG being used as fuel or evaporating (boil-off). For each big pair $p \in \mathcal{P}_v$ we denote by $c_p$ the first contract in the pair, so the volume consistency constraints are: 

\begin{align}
&V_{p}^{start} = V_{p-1}^{end} + V_{c_p} + F_{c_p},& \forall v \in \mathcal{V}, \quad  \forall p \in \mathcal{P}_v: p = (buy, sell) \label{vol_consist1}\\
&V_{p}^{start} = V_{p-1}^{end} - V_{c_p} - W_{c_p} + F_{c_p},& \forall v \in \mathcal{V}, \quad  \forall p \in \mathcal{P}_v: p = (sell, buy) \label{vol_consist2}\\
&V_p^{end} = V_p^{start} - \sum_{c\in C_{t(p)}^v} V_c^v - \sum_{l\in L_{t(p)}^v}lng_l l & \forall v \in \mathcal{V}, \quad  \forall p \in \mathcal{P}_v \label{vol_consist3}
\end{align}
where $t(p)$ is the trip corresponding to the big pair $p$, $C_{t(p)}^v$ are the small discharge contracts that can be visited between the endpoints of $p$.

We ensure the volumes remain in the feasible zone by constraining the partial volumes at each big pair's start and end points.
\begin{align}
    & 0 \leq V_p^{end} \leq V_p^{start} \leq 0.25 vol_v, & \forall v \in \mathcal{V},\quad \forall p \in \mathcal{P}_v: p = (sell, buy)\label{sail_vol_correct1}\\
    & 0.75 vol_v \leq V_p^{end} \leq V_p^{start} \leq vol_v, & \forall v \in \mathcal{V}, \quad \forall p \in \mathcal{P}_v: p = (buy, sell)\label{sail_vol_correct2}
\end{align}
As the LNG volume monotonically decreases between these points, constraints (\ref{sail_vol_correct1}-\ref{sail_vol_correct2}) guarantee the correctness of the sailing volumes on the entire journey.

At the end of the planning horizon, the vessel should have enough LNG to sail to its final destination.
\begin{align}
    &lng_{final}^v \leq V_{p_f}^{end} \leq 0.25 vol_v, & \forall v \in \mathcal{V} \label{final_volume}
\end{align}

The model (\ref{small_discharges_objective}–\ref{final_volume}) is a Mixed Integer Linear Program with binary variables for sailing legs and continuous variables for volumes. Prefixing big contracts constitutes the major difference with the direct Mixed Integer Program formulation from the MIP model presented in \ref{sec:mip}. This trick enables a straightforward way to ensure that the volume onboard remains in the feasible zone during sailing (\ref{sail_vol_correct1}-\ref{sail_vol_correct2}) without any non-linear constraints.

\subsection{Black-box optimization}
\label{sec:black-box}

In order to enhance the performance of the Large Neighborhood Search \textit{(LNS)} consisting of the small discharge insertion, we try out different starting points. Each starting point that we take is a solution of the big-pair model, but with different trip values resulting from different over-delivery penalty $O$ and multi-destination potential reward $R$. 

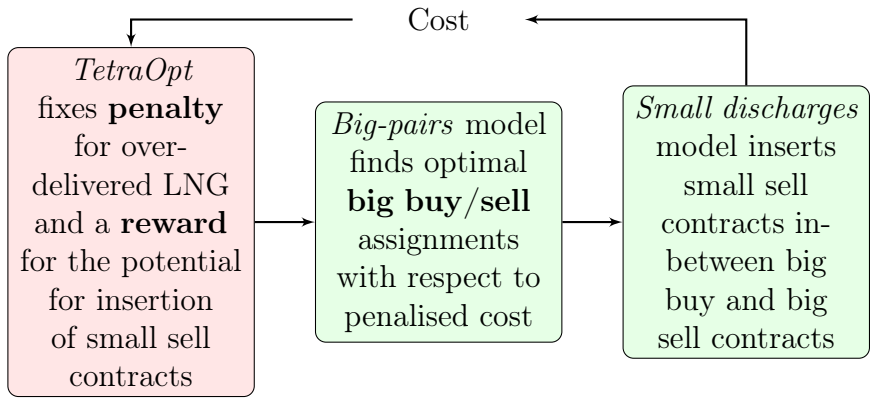
\begin{figure}[htbp]
    \centering
    \begin{tikzpicture}[node distance = 0.8cm]
\node[fill = red!10, text width = 3cm, scheme step] at (0, 0) (penalty){\textit{TetraOpt} fixes \textbf{penalty} \\ for over-delivered LNG \\ and a \textbf{reward} for the potential for insertion of small sell contracts};
\node[fill = green!10, text width = 3cm, scheme step, right = of penalty] (bigpairs){\textit{Big-pairs} model finds optimal \textbf{big buy/sell} assignments with respect to penalised cost};
\node[fill = green!10, text width = 3cm, scheme step, right = of bigpairs] (md) {\textit{Small discharges} model inserts small sell contracts in-between big buy and big sell contracts};
\node[text width = 2cm, above = of bigpairs, align=center, text centered](cost){Cost};
\path[connector] (penalty) -- (bigpairs);
\path[connector] (bigpairs) -- (md);
\path[connector] (md.north) |- (cost);
\path[connector] (cost) -| (penalty.north);
    \end{tikzpicture}
    \caption{The Multistart Large Neighborhood Search heuristic for the LNG shipping problem. The \textit{penalty optimization} (in red) is a black-box optimization problem solved with \textit{TetraOpt}. \textit{Big-pair selection} and \textit{Small discharge insertion} are modeled as ILP and solved with SCIP -- an open-source MILP solver.}
    \label{fig:heuristic_stages}
\end{figure}

Rather than picking the values $O$ and $R$ at random, we consider the problem of optimal selection for $O$ and $R$ as an optimization problem over two continuous variables. These coefficients in the value of the trip (\ref{eq:penalized_cost}) impact the solution in a non-linear and intricate fashion. Therefore, we incline towards a black-box optimization routine to fine-tune these parameters. The objective function corresponds to the best solution in the neighborhood: to get it, we first solve the big-pair model for the selected penalty values, and then we insert small discharges in there. The resulting profit is the cost that we want to optimize. Figure \ref{fig:heuristic_stages} illustrates the whole idea.

\textit{TetraOpt} is a black-box iterative tensor-train solver developed by \textit{Terra Quantum}\footnote{https://terraquantum.swiss/}. It is based on the \textit{maximum volume method} in the construction of the low-rank approximation for multidimensional tensors \citep{TTOpt}, and showed a good performance in protein-protein docking \citep{protein}, conformational sampling of organic molecules \citep{chem_sampling}, neural network architecture optimization \citep{tetraaml,tt_hyp_opt_nn}, and shape optimization \citep{shape}. In our case, each function evaluation corresponds to the optimization of two models and is quite expensive with respect to the runtime requirements. Therefore, we provide the fixed time budget as a parameter for TetraOpt. TetraOpt terminates if the maximum runtime is achieved.

\section{Computational experiments}\label{sec:experiments}
\label{sec:results}

In our experiments, we aim to evaluate i) how profits increase when allowing small discharges, ii) how good LNS is in solving the model with small discharges, and iii) how the guided restarts improve the solution quality. 

\subsection{Data description}

To properly benchmark the quality of our approach, we need data with the so-called \textit{small sell contracts}. We call a sell contract small if it can be served as a small discharge, i.e. if the minimal bound for traded LNG volume is less than $25\%$ of the biggest vessel volume. Due to high sailing costs, it is not economically viable to serve many small contracts in different locations, so the high number of possible small discharges does not necessarily imply that a larger number of them are in the optimal solution. 

Moreover, the small discharges insertion model not only inserts small contracts but also optimizes volumes for big contracts. It means that even if the number of small discharges in data is small (or there are no small discharges at all), our two-stage model may still increase the profit compared to the model where the volumes are predetermined.

Our model uses datasets with i) vessel consumption tables (for example, look at Table \ref{tab:consumption}) ii) predicted free LNG prices for all locations for every day of the planning horizon iii) contract data (type, time windows, and volume bounds) iv) port data v) distances between ports. The planning horizon for the datasets exceeds one year.

\begin{table}[]

\begin{tabular}{l l l l l l}
\hline
Vessel name &
  Fuel mode &
  Is laden &
  \begin{tabular}[c]{@{}l@{}}Speed, \\ knots\end{tabular} &
  \begin{tabular}[c]{@{}l@{}}Consumption, \\ m3/hour\end{tabular} &
  \begin{tabular}[c]{@{}l@{}}Boil-off,\\  m3/hour\end{tabular} \\ \hline
Seal & LNG only & Yes & 20.5 & 9.73E-06 & 8.29E-06 \\
Seal & LNG only & Yes & 20   & 9.06E-06 & 8.29E-06 \\
Seal & LNG only & Yes & 19.5 & 8.52E-06 & 8.29E-06 \\
Seal & LNG only & Yes & 19   & 8.02E-06 & 8.29E-06 \\
Seal & LNG only & Yes & 18.5 & 7.53E-06 & 8.29E-06 \\
Seal & LNG only & Yes & 18   & 7.07E-06 & 8.29E-06 \\
Seal & LNG only & Yes & 17.5 & 6.70E-06 & 6.78E-06 \\
Seal & LNG only & Yes & 17   & 6.34E-06 & 6.78E-06 \\
Seal & LNG only & Yes & 16.5 & 6.01E-06 & 6.78E-06 \\
Seal & LNG only & Yes & 16   & 5.73E-06 & 6.78E-06 \\ \hline
\end{tabular}%

\caption{Example of the consumption table for the imaginary vessel "Seal".}
\label{tab:consumption}
\end{table}

We performed the experiments on the production datasets. 
However, for privacy reasons, we don't report the characteristics of the production instances. Instead, in Table \ref{Art_data}, we describe three artificial instances that we generated with \textit{SDV Python package} \citep{SDV, SDV_software} so they demonstrate the same patterns as the production data. These data instances are approximately the same size but have different contracts, vessels, and prices.

\begin{table}[H]
\setlength{\tabcolsep}{2pt}
\centering
\footnotesize
\caption{Description of the artificial data}
\label{Art_data}
\begin{tabular}{c c c c c c c c c c}
\hline
data & \makecell{\#\\cargoes} & \makecell{\#\\vessels} & \makecell{\#\\buy\\contracts} & \makecell{\#\\sell\\contracts} & \makecell{\#\\small\\contracts} & \multicolumn{2}{c}{\makecell{\#\\flexible\\contracts}} & \makecell{flexible \\ contracts,\\\%}& \makecell{Planning \\ Horizon,\\days} \\
\hline
& & & & & &\makecell{Buy} & \makecell{Sell} &  & \\ 

$A_1$ & 252 & 4 & 52 & 200 & 133 & 52 & 200 & 100 & 885 \\
$A_2$ & 252 & 4 & 52 & 200 & 133  & 52 & 200 & 100 & 885 \\
$A_3$ & 251 & 4 & 51 & 200 & 133 & 51 & 200 & 100 & 885\\
\hline
\end{tabular}
\end{table}

We say that a contract is \textit{flexible} if it is a small contract or if its low volume boundary does not coincide with the upper volume boundary. Such a contract provides the potential for volume optimization or small discharge insertion. In Table \ref{Art_data}, we give the number of flexible contracts. To test the influence of this parameter on the problem solution, we generated datasets \textbf{V2} and \textbf{V4}. Each dataset has a root instance with a high fraction of small sales and diverse volume boundaries. Then, we decrease flexibility by excluding small contracts and equating low volume boundaries to upper ones. All contract prices, locations, and time windows remain the same, as does the fleet. Table \ref{Flexibility_description} describes the resulting datasets.

\begin{table}[H]

\setlength{\tabcolsep}{2pt}
\centering
\footnotesize
\caption{Description of \textbf{V2} and \textbf{V4} datasets}
\label{Flexibility_description}
\begin{tabular}{c c c c c c c c c c}
\hline
data & \makecell{\#\\cargoes} & \makecell{\#\\vessels} & \makecell{\# \\ buy \\ contracts} & \makecell{\#\\sell\\contracts} & \makecell{\#\\small\\contracts} & \multicolumn{2}{c}{\makecell{\#\\flexible\\contracts}} & \makecell{flexible\\contracts,\\ \%}& \makecell{Planning\\Horizon,\\days} \\
\hline
& & & & & &\makecell{Buy} & \makecell{Sell} &  & \\ 

$V4\_026$ & 240 & 4 & 52 & 188 & 36 & 10 & 53 & 26 & 885 \\
$V4\_052$ & 243 & 4 & 52 & 191 & 64  & 24 & 102 & 52 & 885 \\
$V4\_076$ & 248 & 4 & 52 & 196 & 100 & 35 & 154 & 76 & 885\\
$V4\_100$ & 252 & 4 & 52 & 200 & 133 & 52 & 200 & 100 & 885\\
$V2\_027$ & 188 & 2 & 50 & 138 & 22 & 12 & 38 & 27 & 728 \\
$V2\_053$ & 190 & 2 & 50 & 140 & 47  & 24 & 76 & 53 & 728 \\
$V2\_077$ & 194 & 2 & 50 & 144 & 74 & 35 & 115 & 77 & 728\\
$V2\_100$ & 200 & 2 & 50 & 150 & 98 & 50 & 150 & 100 & 728\\
\hline
\end{tabular}
\end{table}

\subsection{Experimental protocol}

We use a n2-standard-8 GCP instance with 8 vCPUs and 32 GB RAM, run by Ubuntu 22.04 for computations. As a MILP solver, we take open-source SCIP \citep{achterberg2009scip}, version 8.0.4. All models are implemented in Python with Pyomo optimization modeling language \citep{hart2011pyomo, bynum2021pyomo}, version 6.6.2. We limit the penalties tuning time for TetraOpt to 12 minutes for all data instances.

\subsection{Results}

We benchmark three approaches. 
In the first approach, we use only the arc-flow big-pairs model without multiple discharges (\ref{cost}-\ref{edge}), similar to the baseline model. 
Secondly, we run the small discharge model, which optimizes volumes and inserts small contracts into the obtained solution, if possible. This corresponds to the large neighborhood search around a single starting point. 
The last approach is our three-component iterative heuristic, where the black-box solver suggests different starting points for the LNS, as described in Section \ref{sec:black-box}. Table \ref{tab:results} shows overall runtimes for each approach, including preprocessing (trips generation) time, which is the same for all three approaches. 
For TetraOpt, we report the highest profit obtained after 12 minutes of penalties tuning. The runtime for TetraOpt in \ref{tab:results} is calculated as a sum of preprocessing time and the actual search time to reach the best found solution (the latter is always less than 12 minutes).

Table \ref{tab:results} shows the runtimes and profits in billions U.S. dollars on the tested data. Given such an order of magnitude, even a slight improvement in profits and losses brings huge benefit.
 
\begin{table}[htbp]
  \centering
  \caption{Computational results for all datasets.}
  \label{tab:results}
  \begin{tabular}{c c c c c c c}
    \hline
    data & \multicolumn{2}{c}{Big pairs} & \multicolumn{2}{c}{Insertion} & \multicolumn{2}{c}{TetraOpt} \\ \hline
    & Profit & Runtime & Profit & Runtime & Profit & Runtime \\ \hline
    $A_1$ & 2.764 & 4 min 26 s & 3.055 & 4 min 47 s & 3.416 & 9 min 14 s \\ 
    $A_2$ & 2.275 & 3 min 14 s & 2.361 & 3 min 15 s & 2.571 & 11 min 56 s \\ 
    $A_3$ & 2.097 & 3 min 2 s & 2.196 & 3 min 4 s & 2.242 & 10 min 18 s \\ 
    $V4\_026$ & 2.397 & 3 min 39 s & 2.805 & 3 min 42 s & 3.136 & 3 min 15 s \\ 
    $V4\_052$ & 2.719 & 4 min 12 s & 3.051 & 4 min 14 s & 3.522 & 12 min 56 s \\ 
    $V4\_076$ & 2.925 & 4 min 47 s & 3.148 & 4 min 49 s & 3.602 & 4 min 50 s \\ 
    $V4\_100$ & 3.071 & 5 min 13 s & 3.338 & 5 min 15 s & 3.589 & 13 min 27 s \\ 
    $V2\_027$ & 1.354 & 1 min 20 s & 1.504 & 1 min 21 s & 1.952 & 5 min 9 s \\ 
    $V2\_053$ & 1.652 & 1 min 26 s & 1.805 & 1 min 27 s & 2.175 & 1 min 40 s \\ 
    $V2\_077$ & 1.829 & 1 min 38 s & 2.019 & 1 min 39 s & 2.279 & 2 min 26 s \\ 
    $V2\_100$ & 1.892 & 1 min 38 s & 2.030 & 1 min 39 s & 2.287 & 2 min 34 s \\ \hline
  \end{tabular}
  \caption*{P\&L in billions of U.S. dollars.}
\end{table}

On production data, we also solved the simplified MIP formulation from \ref{sec:mip}. 
On this dataset, our three-component approach provides a $35\%$ profit improvement over the simple big-pair model in the baseline approach in about six minutes of computational time, while our MIP model with decomposition, solved in 45 minutes with a commercial solver, provides only $27\%$ of profit increase. Both approaches increase profit by allowing multiple discharges; however, our heuristic doesn't require simplifications or time decomposition, keeping more opportunities for cost optimization.\\

We emphasize that, compared to MIP, on production data, our solution approach is not only faster but also provides better solution quality, achieving both simultaneously.
\begin{figure}
    \centering
    \begin{minipage}[b]{0.49\textwidth}
        \includegraphics[width=\textwidth]{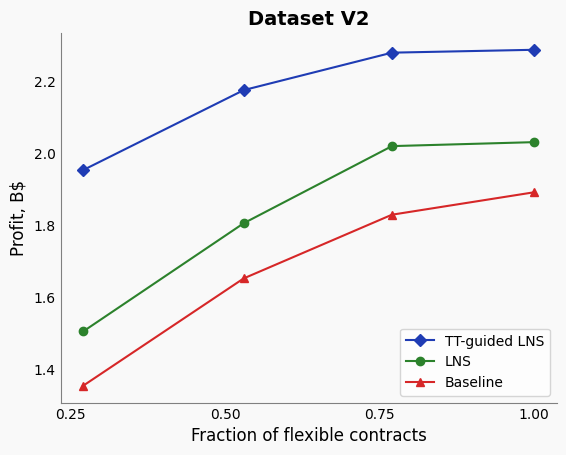}
    \end{minipage}
    \hfill
    \begin{minipage}[b]{0.49\textwidth}
        \includegraphics[width=\textwidth]{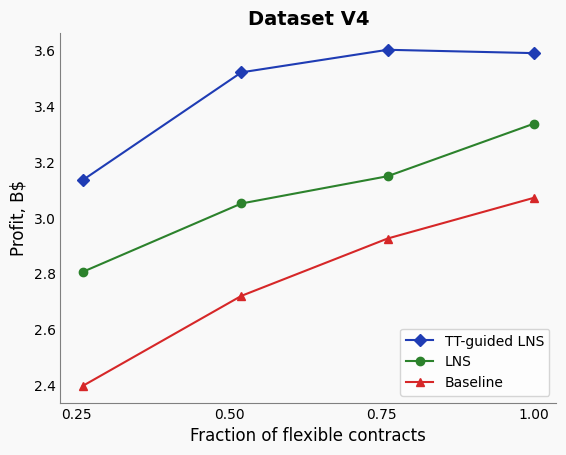}
    \end{minipage}
    \caption{Modeled profit dependency on the data instance for 3 tested approaches. The red line corresponds to the result of the Big-pairs model, the green line -- to one run of the LNS heuristic (taking the baseline as a starting point), and the blue line -- to the best result obtained by Tensor-train guided LNS (with multiple starts).}
    \label{fig:flexibility_exploitation}
\end{figure}

\subsection{Analysis of flexibility exploitation}

Each of the three approaches captures flexible contracts and exploits the trade volume uncertainty to maximize profit. However, our models do this in different ways. The big-pairs selection model works only with predefined feasible trips consisting of one load and one discharge. Thus, the model does not take the full sequence of contracts into consideration and does not optimize the trading volumes through the entire chain of purchases and deliveries. 
Instead, the big-pairs selection model sells the most profitable allowed volume for each chosen laden trip, providing local and, hence, limited optimization for trading volumes. 
Also, the big-pairs approach completely ignores small contracts. 
In contrast, the small discharges insertion model distributes the volume optimally and generates profitable small sales. Nevertheless, it suffers from the suboptimal choice of predefined big contracts.
Our diversified LNS heuristic exploits flexible contracts better than other approaches because it searches for the best big contracts assignment and combines it with effective volume optimization.

Figure (\ref{fig:flexibility_exploitation}) shows the difference in the performance of the three tested approaches. All models benefit from the flexibility of the data; the profit grows in each case. Our heuristic suggests the most effective flexibility utilization in comparison to other tested solutions.

\section{Conclusion}
\label{sec:conclusion}

\begin{figure}
    \centering
\includegraphics[width=\textwidth]{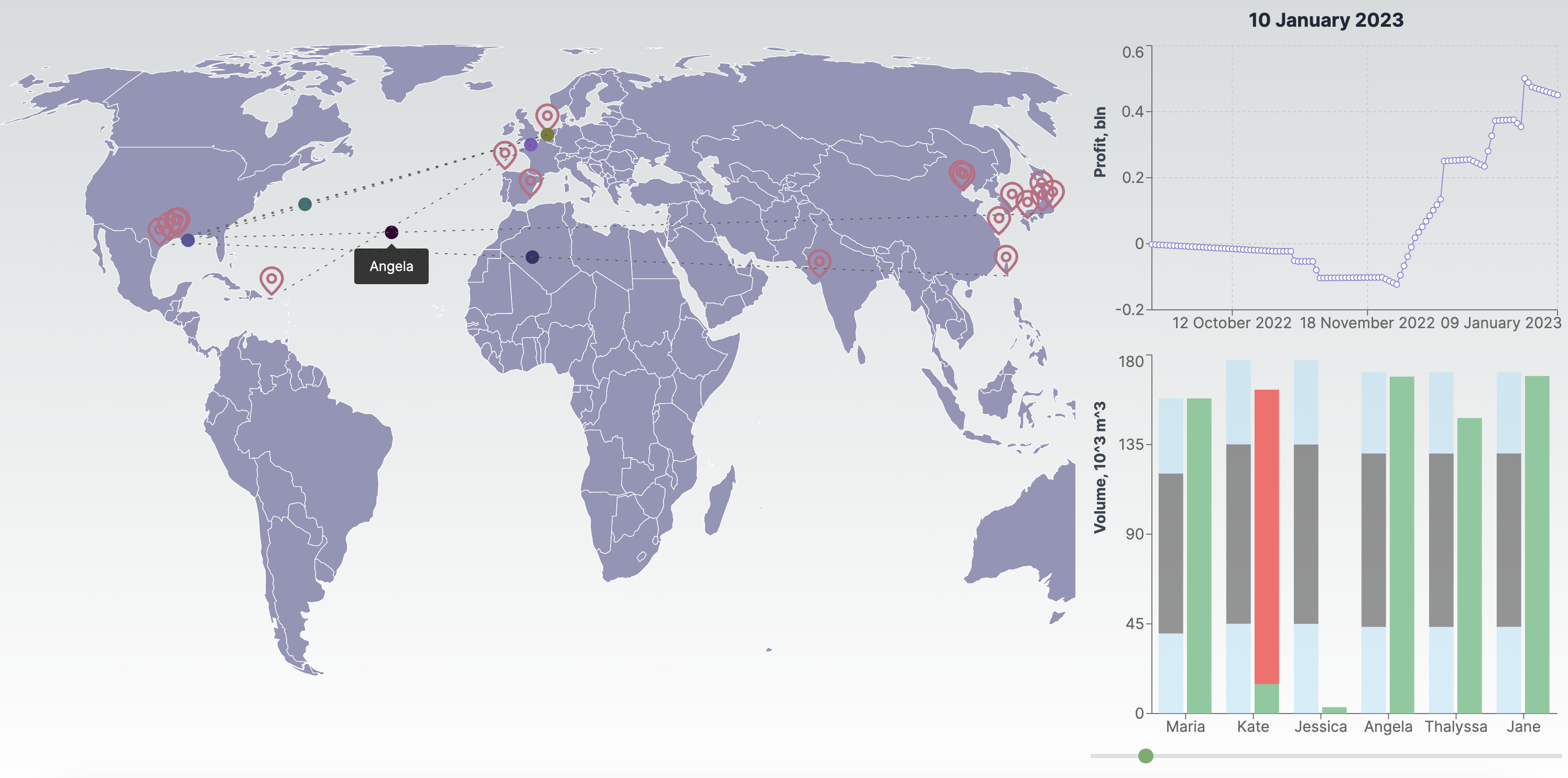}
    \caption{A snapshot of the optimized plan for January 10, 2023. Our visualization tool shows the location of all vessels (in ports or on the journey legs). For each vessel, it tracks the volume onboard; we also plot forbidden volume zones (in gray) and allowed volumes (in blue).}
    \label{fig:demo}
\end{figure}
In this paper, we addressed the new LNG transportation problem, which allows a heterogeneous fleet with restrictions on vessel availability, several types of fuel, and consumption that depends on payload and speed. We considered spot contracts with specified locations and boundaries on time and volume. Our solution takes into account sloshing and multiple discharges. 

We developed a Large Neighborhood Search heuristic, which uses two consecutive MILP models to schedule big trades and small discharges between them for each vessel. The tensor-train optimizer guides the search for the most profitable solution by suggesting new starting points for the LNS by modifying big-pair model parameters. 
Due to the possibility of multiple discharges, our solution outperforms the baseline model without multi-destination by $35\%$ on production data and by $7-44\%$ on artificially generated data with a modest time overhead of several minutes. 

Our algorithm was integrated into a visualization tool showing the snapshots of our solution at each moment in time (see Fig. \ref{fig:demo}). This tool helps to understand the structure of the solution and the role of small discharges, and clearly demonstrates the satisfaction of the volume restrictions. 
Such a tool paves the way towards a full-scale decision-support software for LNG trading and shipping.

Our results demonstrate how to perturb a solution provided by a simple arc-flow model to obtain vessels' paths with a more complex structure and higher profit. 
We expect our heuristic to be helpful for any kind of transportation problem, which may be separated into a simple "core" optimization model and a more complex "perturbing" one.

\section*{Acknowledgments}
We thank Hongda Jiang for fruitful discussions and Jack Lawrence for guiding the project.

\newpage
\bibliographystyle{elsarticle-harv}
\bibliography{refs}

@inproceedings{TTOpt,
 author = {Sozykin, Konstantin and Chertkov, Andrei and Schutski, Roman and Phan, Anh-Huy and CICHOCKI, Andrzej S and Oseledets, Ivan},
 booktitle = {Advances in Neural Information Processing Systems},
 editor = {S. Koyejo and S. Mohamed and A. Agarwal and D. Belgrave and K. Cho and A. Oh},
 pages = {26052--26065},
 publisher = {Curran Associates, Inc.},
 title = {{TTOpt: A Maximum Volume Quantized Tensor Train-based Optimization and its Application to Reinforcement Learning}},
 volume = {35},
 year = {2022}
}

@misc{protein,
  title={Protein-protein docking using a tensor train black-box optimization method},
  author={Morozov, Dmitry and Melnikov, Artem and Shete, Vishal and Perelshtein, Michael},
  primaryClass ={bio-ph},
  archivePrefix={arXiv},
  eprint={2302.03410},
  year={2023}
}

@misc{tetraaml,
      title={Tetra-AML: Automatic Machine Learning via Tensor Networks}, 
      author={A. Naumov and Ar. Melnikov and V. Abronin and F. Oxanichenko and K. Izmailov and M. Pflitsch and A. Melnikov and M. Perelshtein},
      year={2023},
      eprint={2303.16214},
      archivePrefix={arXiv},
      primaryClass={cs.LG}
}

@misc{shape,
      title={Optimization of chemical mixers design via tensor trains and quantum computing}, 
      author={Nikita Belokonev and Artem Melnikov and Maninadh Podapaka and Karan Pinto and Markus Pflitsch and Michael Perelshtein},
      year={2023},
      eprint={2304.12307},
      archivePrefix={arXiv},
      primaryClass={quant-ph}
}

@article{Appelgren1969,
  title = {A Column Generation Algorithm for a Ship Scheduling Problem},
  volume = {3},
  ISSN = {1526-5447},
  DOI = {https://doi.org/10.1287/trsc.3.1.53},
  number = {1},
  journal = {Transportation Science},
  publisher = {Institute for Operations Research and the Management Sciences (INFORMS)},
  author = {Appelgren,  Leif H.},
  year = {1969},
  pages = {53–68}
}

@article{christiansen2013_review,
title = {Ship routing and scheduling in the new millennium},
journal = {European Journal of Operational Research},
volume = {228},
number = {3},
pages = {467-483},
year = {2013},
issn = {0377-2217},
doi = {https://doi.org/10.1016/j.ejor.2012.12.002},
author = {Marielle Christiansen and Kjetil Fagerholt and Bjørn Nygreen and David Ronen}
}

@article{psaraftis_horse,
author = {Psaraftis, Harilaos},
year = {2019},
pages = {1-14},
title = {Ship routing and scheduling: the cart before the horse conjecture},
volume = {21},
journal = {Maritime Economics and Logistics},
doi = {https://doi.org/10.1057/s41278-017-0080-x}
}

@misc{
UNCTAD2023, title={{Review of Maritime Transport 2023}}, 
url={https://unctad.org/publication/review-maritime-transport-2023}, 
author={Hoffmann, Jan and Asariotis, Regina and Assaf, Mark and Bacrot, Celine and Benamara, Hassiba and Hansen, Poul and Hoffmann, Jan and Kulaga, Tomasz and Premti, Anila and Rodríguez, Luisa and Youssef, Frida},
note={{United Nations Publications}, {A}ccessed October 27, 2025}
}

@ARTICLE{harwood_quantum,
  author={Harwood, Stuart and Gambella, Claudio and Trenev, Dimitar and Simonetto, Andrea and Bernal, David and Greenberg, Donny},
  journal={IEEE Transactions on Quantum Engineering}, 
  title={Formulating and Solving Routing Problems on Quantum Computers}, 
  year={2021},
  volume={2},
  number={},
  pages={1-17},
  doi={https://doi.org/10.1109/TQE.2021.3049230}}

@article{li_discretized_speeds,
author = {Li, Lingzi and Ji, Bin and Yu, Samson and Zhou, Saiqi and Fang, Xiaoping},
year = {2022},
pages = {1811},
title = {Branch-And-Price Algorithm for the Tramp Ship Routing and Scheduling Problem Considering Ship Speed and Payload},
volume = {10},
journal = {Journal of Marine Science and Engineering},
doi = {https://doi.org/10.3390/jmse10121811}
}

@article{Homsi2020_UHGS,
title = {Industrial and tramp ship routing problems: Closing the gap for real-scale instances},
journal = {European Journal of Operational Research},
volume = {283},
number = {3},
pages = {972-990},
year = {2020},
issn = {0377-2217},
doi = {https://doi.org/10.1016/j.ejor.2019.11.068},
author = {Gabriel Homsi and Rafael Martinelli and Thibaut Vidal and Kjetil Fagerholt},
}

@article{Rakke_ADP_and_rolling_horizon,
author = {Rakke, Jørgen and Stålhane, Magnus and Moe, Christian and Christiansen, Marielle and Andersson, Henrik and Fagerholt, Kjetil and Norstad, Inge},
year = {2011},
pages = {896-911},
title = {A rolling horizon heuristic for creating a liquefied natural gas annual delivery program},
volume = {19},
journal = {Transportation Research Part C Emerging Technologies},
doi = {https://doi.org/10.1016/j.trc.2010.09.006}
}

@article{LNGScheduler,
author = {Fodstad, Marte and Stremersch, Geert and Hecq, Stéphane and Uggen, Kristin and Rømo, Frode and Lium, Arnt-Gunnar},
year = {2010},
pages = {31-64},
title = {LNGScheduler: A rich model for coordinating vessel routing, inventories and trade in the LNG supply chain},
volume = {3},
journal = {Journal of Energy Markets},
doi = {https://doi.org/10.21314/JEM.2010.049}
}

@InProceedings{Haug_ADP,
author="Haug, Helle V.
and Solum, Sigrid H.
and Warholm, Sanna B.
and Fagerholt, Kjetil
and Li, Mingyu
and Norstad, Inge",
editor="Daduna, Joachim R.
and Liedtke, Gernot
and Shi, Xiaoning
and Vo{\ss}, Stefan",
title={{Planning LNG Annual Delivery Programs with Speed Optimization and Multiple Loading Ports}},
booktitle="Computational Logistics",
year="2023",
publisher="Springer Nature Switzerland",
address="Cham",
pages="170--184",
isbn="978-3-031-43612-3",
doi = {https://doi.org/10.1007/978-3-031-43612-3_10}
}

@techreport{SCIP,
  author = {Ksenia Bestuzheva and Mathieu Besan{\c{c}}on and Wei-Kun Chen and Antonia Chmiela and Tim Donkiewicz and Jasper van Doornmalen and Leon Eifler and Oliver Gaul and Gerald Gamrath and Ambros Gleixner and Leona Gottwald and Christoph Graczyk and Katrin Halbig and Alexander Hoen and Christopher Hojny and Rolf van der Hulst and Thorsten Koch and Marco L{\"u}bbecke and Stephen J. Maher and Frederic Matter and Erik M{\"u}hmer and Benjamin M{\"u}ller and Marc E. Pfetsch and Daniel Rehfeldt and Steffan Schlein and Franziska Schl{\"o}sser and Felipe Serrano and Yuji Shinano and Boro Sofranac and Mark Turner and Stefan Vigerske and Fabian Wegscheider and Philipp Wellner and Dieter Weninger and Jakob Witzig},
  title = {{The SCIP Optimization Suite 8.0}},
  type = {Technical Report},
  institution = {Optimization Online},
  month = {December},
  year = {2021},
  url = {http://www.optimization-online.org/DB\_HTML/2021/12/8728.html},
  note = {{A}ccessed October 27, 2025}
}

@article{Multu_multiple_discharges,
  title={A comprehensive annual delivery program for upstream liquefied natural gas supply chain},
  author={Mutlu, Fatih and Msakni, Mohamed K and Yildiz, Hakan and S{\"o}nmez, Erkut and Pokharel, Shaligram},
  journal={European Journal of Operational Research},
  volume={250},
  number={1},
  pages={120--130},
  year={2016},
  publisher={Elsevier},
  doi = {https://doi.org/10.1016/j.ejor.2015.10.031}
}

@inbook{gronhaug2009_supply_chain,
author = {Grønhaug, Roar and Christiansen, Marielle},
year = {2009},
pages = {195-218},
title = {Supply Chain Optimization for the Liquefied Natural Gas Business},
volume = {619},
isbn = {978-3-540-92943-7},
publisher = {Innovations in Distribution Logistics, Lecture Notes in Economics and Mathematical Systems},
doi = {https://doi.org/10.1007/978-3-540-92944-4_10}
}

@inproceedings{lehr2017comparative,
  title={Comparative threat from LNG and fuel oil maritime accidents},
  author={Lehr, William J and Simecek-Beatty, Debra},
  booktitle={International Oil Spill Conference Proceedings},
  volume={2017},
  number={1},
  pages={3151--3162},
  year={2017},
  doi = {https://doi.org/10.7901/2169-3358-2017.1.3151}
}

@article{coelho2014thirty,
 URL = {http://www.jstor.org/stable/43666993},
 author = {Leandro C. Coelho and Jean-François Cordeau and Gilbert Laporte},
 journal = {Transportation Science},
 number = {1},
 pages = {1--19},
 publisher = {INFORMS},
 title = {Thirty Years of Inventory Routing},
 note = {{A}ccessed October 27, 2025},
 volume = {48},
 year = {2014}
}

@inproceedings{
    SDV,
    title={The Synthetic data vault},
    author={Patki, Neha and Wedge, Roy and Veeramachaneni, Kalyan},
    booktitle={IEEE International Conference on Data Science and Advanced Analytics (DSAA)},
    year={2016},
    pages={399-410},
    doi={https://doi.org/10.1109/DSAA.2016.49},
}

@article{achterberg2009scip,
  title={SCIP: solving constraint integer programs},
  author={Achterberg, Tobias},
  journal={Mathematical Programming Computation},
  volume={1},
  pages={1--41},
  year={2009},
  publisher={Springer},
  doi = {https://doi.org/10.1007/s12532-008-0001-1}
}

@article{zheltkov2020global,
title = {{Global optimization based on TT-decomposition}},
author = {Dmitry Zheltkov and Eugene Tyrtyshnikov},
pages = {247--261},
volume = {35},
number = {4},
journal = {Russian Journal of Numerical Analysis and Mathematical Modelling},
doi = {https://doi.org/10.1515/rnam-2020-0021},
year = {2020},
lastchecked = {2025-10-27}
}

@article{Halvorsen_ADP,
author = {Halvorsen-Weare, Elin and Fagerholt, Kjetil},
year = {2010},
pages = {1-20},
title = {Routing and scheduling in a liquefied natural gas shipping problem with inventory and berth constraints},
volume = {203},
journal = {Annals of Operations Research},
doi = {https://doi.org/10.1007/s10479-010-0794-y}
}

@article{LNS_for_LNG-IRP,
author = {Goel, Vikas and Furman, Kevin and Song, Jin-Hwa and El-Bakry, Amr},
year = {2012},
pages = {},
title = {Large neighborhood search for LNG inventory routing},
volume = {18},
journal = {Journal of Heuristics},
doi = {https://doi.org/10.1007/s10732-012-9206-6}
}

@article{decomposition_LNG-IRP,
author = {Andersson, Henrik and Christiansen, Marielle and Desaulniers, Guy},
year = {2015},
pages = {1-15},
title = {A new decomposition algorithm for a liquefied natural gas inventory routing problem},
volume = {54},
journal = {International Journal of Production Research},
doi = {https://doi.org/10.1080/00207543.2015.1037024}
}

@article{fix-and-relax_heuristic_LNG-IRP,
author = {Uggen, Kristin and Fodstad, Marte and Nørstebø, Vibeke},
year = {2011},
pages = {},
title = {Using and extending fix-and-relax to solve maritime inventory routing problems},
volume = {21},
journal = {Top},
doi = {https://doi.org/10.1007/s11750-011-0174-z}
}

@article{stalhane_adp_heuristic,
author = {Stålhane, Magnus and Rakke, Jørgen and Moe, Christian and Andersson, Henrik and Christiansen, Marielle and Fagerholt, Kjetil},
year = {2012},
pages = {245-255},
title = {A construction and improvement heuristic for a liquefied natural gas inventory routing problem},
volume = {62},
journal = {Computers \& Industrial Engineering},
doi = {https://doi.org/10.1016/j.cie.2011.09.011}
}

@article{sloshing_calcium_carbonate,
author = {Dauzère-Pérès, Stéphane and Nordli, Atle and Olstad, Asmund and Haugen, Kjetil and Koester, Ulrich and Myrstad, Per and Teistklub, Geir and Reistad, Alf},
year = {2007},
pages = {39–51},
title = {Omya Hustadmarmor Optimizes Its Supply Chain for Delivering Calcium Carbonate Slurry to European Paper Manufacturers},
volume = {37},
journal = {Interfaces},
doi = {https://doi.org/10.1287/inte.1060.0276}
}

@article{papaleonidas_decision_support_tool,
  title={An innovative decision support tool for liquefied natural gas supply chain planning},
  author={Papaleonidas, Christos and Lyridis, Dimitrios V and Papakostas, Alexios and Konstantinidis, Dimitris Antonis},
  journal={Maritime Business Review},
  volume={5},
  number={1},
  pages={121--136},
  year={2020},
  publisher={Emerald Publishing Limited},
  doi = {https://doi.org/10.1108/MABR-09-2019-0036}
}

@article{bittante_small_scale,
title = {Optimization of a small-scale LNG supply chain},
journal = {Energy},
volume = {148},
pages = {79-89},
year = {2018},
issn = {0360-5442},
doi = {https://doi.org/10.1016/j.energy.2018.01.120},
author = {A. Bittante and F. Pettersson and H. Saxén}
}

@article{review_2023,
author = {Fagerholt, Kjetil and Hvattum, Lars Magnus and Papageorgiou, Dimitri J. and Urrutia, Sebastián},
title = {Maritime inventory routing: recent trends and future directions},
journal = {International Transactions in Operational Research},
volume = {30},
number = {6},
pages = {3013-3056},
doi = {https://doi.org/10.1111/itor.13313},
year = {2023}
}

@article{qatar_case_study,
title = {An optimization approach to increasing sustainability and enhancing resilience against environmental constraints in LNG supply chains: A Qatar case study},
journal = {Energy Reports},
volume = {8},
pages = {9742-9756},
year = {2022},
issn = {2352-4847},
doi = {https://doi.org/10.1016/j.egyr.2022.07.120},
author = {Sara Al-Haidous and Rajesh Govindan and Adel Elomri and Tareq Al-Ansari},
}

@article{msakni_short_term,
title = {Short-term planning of liquefied natural gas deliveries},
journal = {Transportation Research Part C: Emerging Technologies},
volume = {90},
pages = {393-410},
year = {2018},
issn = {0968-090X},
doi = {https://doi.org/10.1016/j.trc.2018.03.013},
author = {Mohamed Kais Msakni and Mohamed Haouari},
}

@article{rakke_delivery_patterns,
  title={A new formulation based on customer delivery patterns for a maritime inventory routing problem},
  author={Rakke, J{\o}rgen Glomvik and Andersson, Henrik and Christiansen, Marielle and Desaulniers, Guy},
  journal={Transportation Science},
  volume={49},
  number={2},
  pages={384--401},
  year={2015},
  publisher={INFORMS},
  doi = {https://doi.org/10.1287/trsc.2013.0503}
}

@article{gronhaug_bnp,
author = {Grønhaug, Roar and Christiansen, Marielle and Desaulniers, Guy and Desrosiers, Jacques},
year = {2010},
pages = {400-415},
title = {A Branch-and-Price Method for a Liquefied Natural Gas Inventory Routing Problem},
volume = {44},
journal = {Transportation Science},
doi = {https://doi.org/10.2307/25769507}
}

@misc{QAOA_original,
  title={A quantum approximate optimization algorithm},
  author={Farhi, Edward and Goldstone, Jeffrey and Gutmann, Sam},
  year={2014},
  eprint={1411.4028},
  archivePrefix={arXiv},
  primaryClass={quant-ph}
}

@misc{farhi_qa,
  title={Quantum computation by adiabatic evolution},
  author={Farhi, Edward and Goldstone, Jeffrey and Gutmann, Sam and Sipser, Michael},
  eprint={0001106},
  archivePrefix={arXiv},
  primaryClass={quant-ph},
  year={2000}
}

@ARTICLE{Dwave_issues,
AUTHOR={Bian, Zhengbing and Chudak, Fabian and Israel, Robert and Lackey, Brad and Macready, William G. and Roy, Aidan},   
TITLE={Discrete optimization using quantum annealing on sparse Ising models},      
JOURNAL={Frontiers in Physics},      
VOLUME={2},           
YEAR={2014},              
DOI={https://doi.org/10.3389/fphy.2014.00056},      
ISSN={2296-424X}
}

@book{bynum2021pyomo, title={Pyomo--optimization modeling in python}, author={Bynum, Michael L. and Hackebeil, Gabriel A. and Hart, William E. and Laird, Carl D. and Nicholson, Bethany L. and Siirola, John D. and Watson, Jean-Paul and Woodruff, David L.}, edition={Third}, volume={67}, year={2021}, publisher={Springer Science \& Business Media} }

@article{hart2011pyomo, title={Pyomo: modeling and solving mathematical programs in Python}, author={Hart, William E and Watson, Jean-Paul and Woodruff, David L}, journal={Mathematical Programming Computation}, volume={3}, number={3}, pages={219--260}, year={2011}, publisher={Springer}, doi={10.1007/s12532-011-0026-8}}

@misc{SCIP_Optimization_Suite,
  author       = {{SCIP Optimization Suite Development Team}},
  title        = {SCIP Optimization Suite (Version 8.0.4) [Computer software]},
  year         = {2023},
  note    = {Zuse Institute Berlin (ZIB)},
  howpublished = {\url{https://scipopt.org/}},

}

@misc{SDV_software,
  author       = {{SDV Development Team}},
  title        = {Synthetic Data Vault (SDV) (Version 1.4.0) [Computer software]},
  year         = {2023},
  note = {DataCebo, Inc.},
  howpublished         = {\url{https://github.com/sdv-dev/SDV}},
}

@article{chem_sampling,
author = {Zurek, Christopher and Mallaev, Ruslan A. and Paul, Alexander C. and van Staalduinen, Nils and Pracht, Philipp and Ellerbrock, Roman and Bannwarth, Christoph},
title = {{Tensor Train Optimization for Conformational Sampling of Organic Molecules}},
journal = {{Journal of Chemical Theory and Computation}},
volume = {21},
number = {3},
pages = {1459-1475},
year = {2025},
doi = {https://doi.org/10.1021/acs.jctc.4c01275},
}

@Article{tt_hyp_opt_nn,
AUTHOR = {Naumov, A. and Melnikov, A. and Perelshtein, M. and Melnikov, Ar. and Abronin, V. and Oksanichenko, F.},
TITLE = {{Tensor Network Methods for Hyperparameter Optimization and Compression of Convolutional Neural Networks}},
JOURNAL = {{Applied Sciences}},
VOLUME = {15},
YEAR = {2025},
NUMBER = {4},
ARTICLE-NUMBER = {1852},
ISSN = {2076-3417},
DOI = {https://doi.org/10.3390/app15041852}
}

\newpage

\appendix

\section{Nodes-based optimization model}
\label{sec:mip}

\subsection{Assumptions}
\label{extra-assumptions}

Initially, to build the draft of the problem solution, we made several more assumptions:

\begin{enumerate}
    \item The only type of fuel available is LNG onboard itself
    \item Each time window for the contract includes an integer number of days. Each vessel arrives in the middle of the contract time window (rounded down if time window duration in days is odd)
    \item The Optimization period is considerably smaller than the rent period for each vessel.
    \item Boil-off gas rate is independent of vessel speed or load.
\end{enumerate}

Assumptions 1 and 4 started from the fact that burned-out LNG volume is usually small in comparison to trading volumes. That is why we calculated the LNG waste during the voyage approximately. We accepted assumption 2, because all contracts in the data indeed have integer time window durations, and most of them have time windows that last 1 day. Thus, the moment when the vessel arrives to complete the contract is small in comparison to travel times between different ports. Due to assumption 3, we do not consider the final locations of the vessels: the vessels always have enough time after the last contract to reach the desired port.
Then, we built the nodes-based model, which solves the problem in one stage, and developed a simple decomposition technique for it.

\subsection{Full-data model}
\label{full-data}

In this section, we describe a MIP model without decomposition. However, it turned out to be impractical because our machine did not have enough memory to process the solution for production data, since the data had a planning horizon of several years. We describe how we decompose the model to make half-year subproblems computationally tractable in \ref{decomposition}.

\subsubsection{Preprocessing and Sets}

As an input to our model, we had the following sets:

\begin{itemize}
\item There is a set of discharge contracts $I$ and set of loading contracts $J$. The union of these sets yields all possible contracts $\mathfrak{I}$, regardless of the contract type.  Each contract $i\in\mathfrak{I}$ has a price value $P_i$. We also have a set of slots $\mathfrak{S} \subset I.$ -- the sell contracts with required volume lower bound equal to $0$. 
    \item $V_i^{min}$ and $V_i^{max}$ -- the minimum and the maximum demands of the $i$-th contract respectively. If the contract is a sell contract, then both bounds are not positive; otherwise, they are not negative. We have $V_i^{max}=0$, if and only if $i \in \mathfrak{S}.$ 
\item $S_{ij}$ is the distance between a pair of contracts $i$, $j$.
\item $v_k^{max}$ -- the maximum possible speed for each vessel $k$, used to estimate the feasibility of the  trip.
\item $\mathtt{V}_k^{max}$ -- the maximum capacity of the $k$-th vessel. 
\item $R_k^{LNG}(v_k)$ -- the daily LNG consumption rate for $k$-th vessel as a given function of vessel speed $v$. The vessel speed is a discrete value; thus, we have $v_k\in\mathfrak{Y}_k$, where $\mathfrak{Y}_k$ is the set of available speed modes for the $k$-th vessel. Note that there is a minimum daily LNG consumption rate $R_k^{LNG}(0)$ even if a vessel is not moving. This is the so-called gas boiling rate caused by gas evaporation.
\item $r_{i}$ and $d_{i}$ are the release date and deadline of each contract $i$, given in days from the beginning of optimization interval. 

\end{itemize}

Also, we create other sets:

\begin{itemize}
    
\item The moving time between contracts $i$ and $j$, $T_{ij}$ (in days), should be calculated in the following way:
\begin{equation}
    T_{ij} = r_j - r_i + \lfloor\frac{d_j - r_j}{2}\rfloor - \lfloor\frac{d_i - r_i}{2}\rfloor - 1,
\end{equation}
where $\lfloor\cdot\rfloor$ is a rounding down function. $r_j$, $d_j$, $r_i$, $d_i$ are not necessary integers. This time corresponds to the moving time between the middle of time windows, and one day is subtracted for LNG loading/discharge. Remember that according to our assumptions, $d_j-r_j$ and $d_i-r_i$ are integers \ref{extra-assumptions}. $T_{ij}$ may be negative, so we do not consider such a pair of contracts as feasible.

\item Thus, we have a fixed time $T_{ij}$ between different pairs of contracts $(i, j)$ and the fixed (minimally required) average speed; we can calculate the consumption rate $R_{ijk}^{LNG}$ as a function of this average speed:

\begin{equation}
    \widetilde{v_{ij}} = \frac{S_{ij}}{T_{ij}},
\end{equation}
where $\widetilde{v_{ij}}$ is the continuous value of the speed. But we should choose the speed from a discrete set $\mathfrak{Y}_k$ for each vessel $k$. We define this value as the closest available value, greater than $\widetilde{v_{ij}}$:
\begin{equation}
    v_{ijk} = \min\{v_k\in\mathfrak{Y}_k|v_k\geq \widetilde{v_{ij}}\},
\end{equation}
if $0\leq \widetilde{v_{ij}}\leq v_k^{max}$. Otherwise, we can assign to the $v_{ijk}$ any special value (infeasibility indicator) corresponding to the situation when $\widetilde{v_{ij}}$ is negative (moving in the past), or it is greater than the maximum possible speed. For example, we define $v_{ijk}=0$ in such cases. After that, we calculate the consumption from the given table:
\begin{equation}
    R_{ijk}^{LNG} = R_k^{LNG}(v_{ijk}),
\end{equation}
if $v_{ijk} \neq 0$. In the opposite case, consumption is unnecessary because the set of indices $ijk$ is impossible.

For each vessel $k$, we can introduce the set of impossible contracts $\mathfrak{O}_k$, which can not be completed by this vessel because of the limitation of the vessel's capacity:

\begin{multline}
    \mathfrak{O}_k = \Bigl\{i\Bigl | \left(V_i^{min} > 0.985\mathtt{V}_k^{max}, i\in J\right) \\ \text{or}\quad\left(\left|V_i^{max}\right| > 0.985\mathtt{V}_k^{max}, i \in I\right)\Bigl\},
\end{multline}

and also we can introduce another set of contracts $\mathfrak{M}_{ki}$, which can not be visited right after cargo $i$ by the vessel $k$. 
If $i, j\in J$, then we should include all contracts satisfying the constraint:

\begin{equation}
    V_i^{min} + V_j^{min} - R_{ijk}T_{ij} > 0.985\mathtt{V}_k^{max}.
\end{equation}
Similarly, if $i, j \in I$, then $j\in\mathfrak{M}_{ki}$ if:

\begin{equation}
    \left|V_i^{max}\right|+\left|V_j^{max}\right| + R_{ijk}T{ij} > 0.985\mathtt{V}_k^{max}.
\end{equation}

Finally, if $R_{ijk}T_{ij} > 0.25\mathtt{V}_k^{max}$, then regardless of the type of the contract, we should include contract $j$ in the set because otherwise, we will violate the constraint on admissible filling percentage for the vessel.

If the vessel $k$ successfully arrived to load (or discharge) the contract $i$ at the moment $q$, then it can not move to the contract $j$ in time step $q+1$, if $j\in\mathfrak{M}_{ki}$. 

\item We introduce the set $\mathfrak{A}_{ki}$ of all contracts $j$, which are impossible for visiting by the vessel $k$ in any time step later than contract $i$ was visited. 
We form it from all the contracts $j$, for which $v_{ijk} = 0$ (infeasibility indicator), regardless of the type of contract.

\item We will denote the union of the sets $\mathfrak{A}_{ki}$ and $\mathfrak{M}_{ki}$ as $\mathfrak{C}_{ki} = \mathfrak{A}_{ki} \cup \mathfrak{M}_{ki}.$

\item We introduce the ancillary sets $\mathfrak{L}_k$ of all contracts, which can not be the first visiting point for the vessel $k$. 
 
We define it with all loading contracts $i\in J$ if:

\begin{equation}
    V_i^{min} > 0.985\mathtt{V}_k^{max},
\end{equation}

and all discharge contracts $i\in I$.
\end{itemize}

\subsubsection{Variables}

\begin{itemize}
    \item
    \begin{equation}
        \mathbf{x}_{ikq} = 
        \begin{cases}
            1, & \text{if the $k$-th vessel processed $i$-th contract in the time step $q$,}\\
            0, & \text{otherwise}
        \end{cases}
    \end{equation}
    \item $\mathbf{V}_{ikq}$ -- the volume of LNG to be manipulated with by the $k$-th vessel at the $q$-th moment and $i$-th contract. It is a continuous variable within the range
    \begin{equation}
    \begin{gathered}
        V_i^{min} \leq \mathbf{V}_{ikq} \leq 0 \quad \quad \forall i \in I,\\
        0 \leq \mathbf{V}_{ikq} \leq V_i^{max} \quad \quad \forall i \in J.
    \end{gathered}
    \end{equation}
    
    Depending on the contract type, it can be negative (if we sell LNG) or positive (if we buy LNG).
    
    \item We also have a slack variable $\bm{\kappa}_{kq}$, which is equal to $1$, if the loaded volume of LNG is more than $75\%$ of a maximum value for the $k$-th vessel on a $q$-th step, but less than $98.5\%$, and $\bm{\kappa}_{kq}=0$ if the loaded volume of LNG is less than $25\%$ (another loaded volume is not considered):
    \begin{equation}
        \bm{\kappa}_{kq} = 
        \begin{cases}
            1 & \text{if} \quad 0.75 \mathtt{V}_k^{max} \leq L_{kq} \leq 0.985\mathtt{V}_k^{max}\\
            0 & \text{if} \quad 0 \leq L_{kq} \leq 0.25 \mathtt{V}_k^{max},
        \end{cases}
    \end{equation}
    where $L_{kq}$ is a loaded volume of LNG on a $k$-th vessel at $q$-th moment.
\end{itemize}

\subsubsection{MIP}

\paragraph{Constraint on the volume of LNG}

We must restrict the trading volume because of the vessel's limited capacity. It yields:

\begin{equation}
   0 \leq L_{kq} + \sum_{i\in \mathfrak{I}}\mathbf{V}_{ikq} \leq \mathtt{V}_k^{max} \quad \forall k, q.
\end{equation}
where $\mathtt{V}_k^{max}$ is the maximum available trading volume on the $k$-th vessel, which can be fulfilled by LNG and $L_{kq}$ is the volume of LNG on the k-th vessel before q-th moment.

The current volume of LNG on the $k$-th vessel before $q$-th moment can be calculated as the sum of all trade actions before the $q$-th moment and all wastes while moving from one terminal to another:

\begin{equation}
    L_{kq} = \sum_{i, q'<q}\mathbf{V}_{ikq'}
    -\sum_{i, q'<q}\sum_{j \not\in \mathfrak{C}_{ki} }\left(R_{ijk}^{LNG}T_{ij} + R_k^{LNG}(0)\right)\mathbf{x}_{ikq'} \mathbf{x}_{jk,q'+1},
\end{equation}

where $R_{ijk}^{LNG}$ is the consumption rate of the vessel $k$ on a road between $i$-th and $j$-th contracts, $T_{ij}$ is the time of movement. We included the boil of gas term $R_k^{LNG}(0)$ in the sum because of the boil of gas during the loading or discharge.

One extra constraint says that if the vessel $k$ doesn't visit the contract $i$ at the moment $q$, then the corresponding volume variable is equal to $0$:

\begin{equation}
\begin{aligned}
    V_i^{min}\mathbf{x}_{ikq} &\leq \mathbf{V}_{ikq} \leq V_i^{max}\mathbf{x}_{ikq}, \quad \quad i \not \in \mathfrak{S},\\
    V_i^{min}\mathbf{x}_{ikq} &\leq \mathbf{V}_{ikq} \leq -R_k^{LNG}(0)\mathbf{x}_{ikq}, \quad \quad i \in \mathfrak{S}.\\
\end{aligned}
\end{equation}

We separate slots and other contracts because it is practically unreasonable to sell any small amount of LNG. We postulate above that we sell at least a comparable amount of LNG to the amount that can boil out while we discharge it.

Finally, the LNG volume on the $k$-th vessel before $q$-th moment (we denote it as $L_{kq}$ should satisfy the bounds:

\begin{equation}
    0.75 \mathtt{V}_k^{max}\bm{\kappa}_{kq} \leq L_{kq} \leq 0.985 \mathtt{V}_k^{max}\bm{\kappa}_{kq} + 0.25\mathtt{V}_{k}^{max}(1 - \bm{\kappa}_{kq}).
\end{equation}

\paragraph{Cost}

We can conclude, using all the above, that the final cost should be expressed as:

\begin{equation}
\begin{gathered}
    C = -\sum_{i,k,q}P_i\mathbf{V}_{ikq}  \to \max.
\end{gathered}
\end{equation}

\paragraph{Physical constraints of different voyages}

If the vessel $k$ successfully arrived to load (or discharge) the contract $i$ at the moment $q$, then it should not move to the contract $j$ at the moment $q+1$, if $j\in\mathfrak{M}_{ki}$. We express this statement as the following constraint:

\begin{equation}
    \mathbf{x}_{ikq} + \mathbf{x}_{jk,q+1} \leq 1 \quad \forall i, q, k, j \in \mathfrak{M}_{ki}.
\end{equation}

If the vessel $k$ successfully arrived to load (or discharge) the contract $i$ at the moment $q$, then it should not move to the contract $j$ at the moment $q'>q$, if $j\in\mathfrak{A}_{ki}$:

\begin{equation}
    \mathbf{x}_{ikq} + \mathbf{x}_{jk,q'} \leq 1 \quad \forall i, q, k, q'>q, j \in \mathfrak{A}_{ki}.
\end{equation}

Every particular contract should not be taken more than once.

\begin{equation}\label{not_more_than_once}
    \sum_{k, q}\mathbf{x}_{ikq} \leq 1 \quad \forall i.
\end{equation}

At each moment, each vessel should not take more than one contract.

\begin{equation}
    \sum_{i}\mathbf{x}_{ikq} \leq 1 \quad \forall k, q.
\end{equation}

If we don't take any contract, then we stop:

\begin{equation}
    \sum_{i}\mathbf{x}_{ik, q+1} \leq \sum_{i}\mathbf{x}_{ikq} \quad \forall k, q.
\end{equation}

Finally, we need to take into account impossible sets $\mathfrak{O}_k$ and $\mathfrak{L}_k$, by their definition:

\begin{equation}
    x_{ikq} = 0 \quad \forall k, \forall q, \forall i \in \mathfrak{O}_k,
\end{equation}

\begin{equation}
    x_{ik0} = 0 \quad \forall k, \forall i \in \mathfrak{L}_k.
\end{equation}

The last two constraints actually reduce the number of variables in the problem.

\paragraph{Model summary}

\begin{itemize}
    \item Cost function:
    
    \begin{equation}\tag{$C_0$}
    \begin{gathered}
    C = -\sum_{i,k,q}P_i\mathbf{V}_{ikq}  \to \max.
    \end{gathered}
    \end{equation}

    \item Constraints on volume of LNG:
    
    \begin{equation}\tag{$C_{01}$}
    \begin{gathered}
        V_i^{min} \leq \mathbf{V}_{ikq} \leq 0 \quad \quad \forall i \in I,\\
        0 \leq \mathbf{V}_{ikq} \leq V_i^{max} \quad \quad \forall i \in J.
    \end{gathered}
    \end{equation}
    
    \begin{equation}\tag{$C_1$}
    L_{kq} = \sum_{i, q'<q}\mathbf{V}_{ikq'}
    -\sum_{i, q'<q}\sum_{j \not\in \mathfrak{C}_{ki} }\left(R_{ijk}^{LNG}T_{ij} + R_k^{LNG}(0)\right)\mathbf{x}_{ikq'} \mathbf{x}_{jk,q'+1}
    \end{equation}
    
    \begin{equation}\tag{$C_2$}
     0 \leq L_{kq} + \sum_{i\in \mathfrak{I}}\mathbf{V}_{ikq} \leq \mathtt{V}_k^{max}.
    \end{equation}
    
    \begin{equation}\tag{$C_3$}
    \begin{aligned}
    V_i^{min}\mathbf{x}_{ikq} &\leq \mathbf{V}_{ikq} \leq V_i^{max}\mathbf{x}_{ikq}, \quad \quad i \not \in \mathfrak{S},\\
    V_i^{min}\mathbf{x}_{ikq} &\leq \mathbf{V}_{ikq} \leq -R_k^{LNG}(0)\mathbf{x}_{ikq}, \quad \quad i \in \mathfrak{S}.\\
    \end{aligned}
    \end{equation}
    
    \begin{equation}\tag{$C_4$}
    0.75 \mathtt{V}_k^{max}\bm{\kappa}_{kq} \leq L_{kq} \leq 0.985 \mathtt{V}_k^{max}\bm{\kappa}_{kq} + 0.25\mathtt{V}_{k}^{max}(1 - \bm{\kappa}_{kq}).
    \end{equation}
    
    \item Constraints on the possibility of different voyages:
    
    \begin{equation}\tag{$C_5$}
    \begin{gathered}
    \mathbf{x}_{ikq} + \mathbf{x}_{jk,q+1} \leq 1 \quad \forall k, i, q, j \in \mathfrak{M}_{ki},\\
    \mathbf{x}_{ikq} + \mathbf{x}_{jkq'} \leq 1 \quad \forall k, i, q, j \in \mathfrak{A}_{ki}, q' > q.
    \end{gathered}
    \end{equation}
    
    \begin{equation}\tag{$C_7$}
    \begin{aligned}
        \sum_{k, q}\mathbf{x}_{ikq} &\leq 1 \quad \forall i,\\
        \sum_{i}\mathbf{x}_{ikq} &\leq 1 \quad \forall k, q, \\
    \end{aligned}
    \end{equation}
    
    \begin{equation}\tag{$C_8$}
        \sum_{i}\mathbf{x}_{ik, q+1} \leq \sum_{i}\mathbf{x}_{ikq} \quad \forall k, q.
    \end{equation}
    
    \begin{equation}\tag{$C_9$}
    x_{ikq} = 0 \quad \forall k, \forall q, \forall i \in \mathfrak{O}_k,
\end{equation}

\begin{equation}\tag{$C_{10}$}
    x_{ik0} = 0 \quad \forall k, \forall i \in \mathfrak{L}_k.
\end{equation}
\end{itemize}

\subsection{Decomposition for nodes-based model}
\label{decomposition}

\begin{figure}[H]
    \centering
    \includegraphics[width=1.0\textwidth]{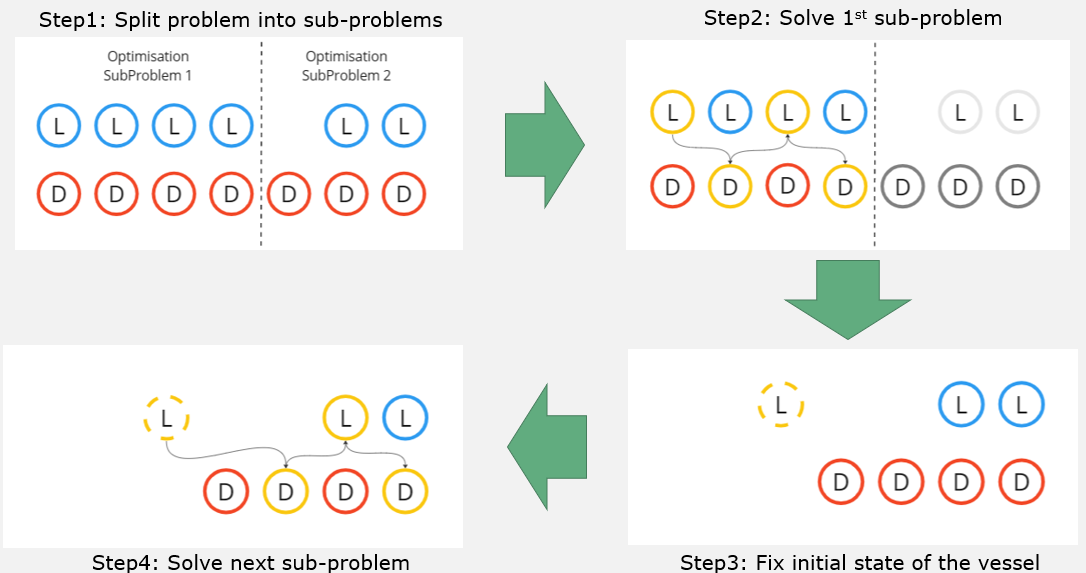}
    \caption{Illustration of the decomposition technique for one vessel. Loading contracts are blue; discharge contracts are red. We mark the contracts taken in the solution with yellow and the contracts from the next subproblem with gray. Note that the vessel may not take the last discharge contracts in the previous subproblem after some more profitable contracts appear in the next subproblem, because we always use the loading contract as an initial one.}
    \label{fig:decomposition}
\end{figure}

We implement a decomposition technique, splitting the time horizon into 6-month periods. Figure \ref{fig:decomposition} explains the idea. We solve the 1st subproblem with a node-based model, then we fix the last loading location for all the vessels and use it as a starting point for the next subproblem. All initial locations and volumes are calculated straightforwardly.

We use a loading contract as a start contract for the vessel, because then the vessel will have enough fuel to move. If it starts in discharge, it may turn out that it can not take any new contract from the second half of the year, and the entire chain of contracts breaks.

The decomposed model requires much less memory and brings more profit than big pairs selection; however, it is not solved optimally in several hours. That is why we decided to implement another solution described in the paper's main part.

\end{document}